\documentclass[10pt]{article}
\usepackage{array}
\usepackage[latin1]{inputenc}

\usepackage{amssymb}
\usepackage{amsfonts}
\newcommand{\bE}{{\mathbb{E}}}

\newcommand{\bJ}{{\mathbb{J}}}

\newcommand{\cC}{{\cal C}}
\newcommand{\cG}{{\cal G}}
\newcommand{\cI}{{\cal I}}
\newcommand{\cF}{{\cal F}}
\newcommand{\cB}{{\cal B}}

\newcommand{\cM}{{\cal M}}

{

\newcommand{\cH}{{\cal H}}

\setlength{\parskip}{10pt plus 2pt minus 1pt}
\setlength{\parindent}{0pt}

\setcounter{section}{-1}

\begin{document}
\begin{center}
\textbf{\LARGE Orthocentric simplices and their centers}
\end{center}

\begin{center}
Allan L. Edmonds$^a$, Mowaffaq Hajja$^{b1}$, Horst
Martini$^c$\footnote{The second named author was supported by a
grant from Yarmouk University, and the third named author by a DFG
grant.}

$^a$ Department of Mathematics\\Indiana
University\\Bloomington, IN 47405\\USA\\[0.3cm]
$^b$ Department of Mathematics\\Yarmouk
University\\Irbid\\JORDAN\\[0.3cm]
$^c$ Faculty of Mathematics\\Chemnitz University of
Technology\\09107 Chemnitz\\GERMANY
\end{center}

\begin{abstract}
A simplex is said to be orthocentric if its altitudes intersect in
a common point, called its orthocenter. In this paper it is proved
that if any two of the traditional centers of an orthocentric
simplex (in any dimension) coincide, then the simplex is regular.
Along the way orthocentric simplices in which all facets have the
same circumradius are characterized, and the possible barycentric
coordinates of the orthocenter are described precisely. In
particular these barycentric coordinates are used to parametrize
the shapes of orthocentric simplices. The substantial, but
widespread, literature on orthocentric simplices is briefly
surveyed in order to place the new results in their proper
context, and some of the previously known results are given
new proofs from the present perspective.

\vspace{0.3cm}

\end{abstract}

\section{Introduction}

This paper is a study of the geometric consequences of assumed
coincidences of \emph{centers} of a \emph{$d$-dimensional
orthocentric simplex} (or, simply, \emph{orthocentric
$d$-simplex}) $S$ in the $d$-dimensional Euclidean space, $d \ge
3$, i.e., of a $d$-simplex $S$ whose $d+1$ altitudes have a common
point $\cH$, called the \emph{orthocenter} of $S$. The centers
under discussion are the \emph{centroid} $\cG$, the
\emph{circumcenter} $\cC$ and the \emph{incenter} $\cI$ of $S$.
For triangles, these centers are mentioned in Euclid's
\emph{Elements}, and in fact they are the only centers mentioned
there. It is interesting that the triangle's \emph{orthocenter}
$\cH$, defined as the intersection of the three altitudes, is
never mentioned in the \emph{Elements}, and that nothing shows
Euclid's awareness of the fact that the three altitudes are
concurrent, see \cite{Gar}, p. 58. It is also worth mentioning
that one of the most elegant proofs of that concurrence is due to
C. F. Gauss, and A. Einstein is said to have prized this
concurrence for its nontriviality and beauty. However, in contrast
to the planar situation, the $d+1$ altitudes of a $d$-simplex are
not necessarily concurrent if $d \ge 3$. We may think of this as a
first manifestation of the reality that general $d$-simplices, $d
\ge 3$, \emph{do not have all the nice properties that triangles
have}.

It is natural that, besides $\cG, \cC$, and $\cI$, we will also
consider the orthocenter $\cH$ of a $d$-simplex $S$ regarding its
coincidence with the other three centers. It is well-known that
for $d = 2$ the coincidence of any two of the four mentioned
centers yields a regular (or equilateral) triangle; see
\cite{Liu}, page 78, and for triangle centers in general we refer
to \cite{B-M} and \cite{Ki}. For $d \ge 3$ this is no longer true,
i.e., only weaker degrees of regularity are obtained, see
\cite{Ed} and \cite{E-H-M} for recent results on this. One of
these weaker degrees is the \emph{equifacetality} of a $d$-simplex
$S$, i.e., the congruence of its $(d-1)$-faces, which does not
imply regularity, see \cite{M-W} and, for a deeper study of
equifacetal simplices, \cite{Ed}. From \cite{E-H-M}, Theorem 3.2.,
it follows that equifacetality implies $\cG = \cC = \cI$, and
there it is also shown that the opposite implication, although
true for $d \in \{2,3\}$, does not hold for $d=4$ (and expectedly
also not for $d \ge 5$). Also the coincidence of two of these
three centers does not imply that all three coincide.
Specifically, the existence of a non-equifacetal $4$-simplex with
$\cG = \cC = \cI$ follows from \cite{Ha-Wa}, \cite{E-H-M}, Theorem
3.2 (iv) and Theorem 4.5, while the existence of $4$-simplices
with $\cG = \cC \not= \cI\,, \, \cC = \cI \not= \cG$, and $\cG =
\cI \not= \cC$ follows from \cite{Ha-Wa}, \cite{E-H-M}, Theorems
3.4 and 4.3, respectively. This emphasizes once again the feeling
that arbitrary $d$-simplices, $d \ge 3$, are not the most faithful
generalizations of triangles. It turns out that orthocentric
$d$-simplices do resemble triangles closely in significant ways.
Namely, we will prove that if $S$ is orthocentric, then
\[
\begin{array}{lll}
S \mbox{ is regular } & \Longleftrightarrow & S \mbox{ is equifacetal}\,,\\
& \Longleftrightarrow & \cG = \cC = \cI = \cH\,,\\
& \Longleftrightarrow & \mbox{any two of these four centers coincide}.
\end{array}
\]
Further similar results, related to other degrees of regularity of $S$, will be added, and also for
a class of special orthocentric simplices, called
\emph{rectangular simplices} (for which $\cH$ is a vertex), various results will be presented.

The class of orthocentric simplices has a long history. The
literature on these special polytopes is large, and  there is no
satisfactory survey showing the current state of knowledge about
them. We therefore give below a short survey of the literature on
orthocentric simplices. The reader will observe that there are
still many possible ways to find new properties of this
interesting class of simplices.

\section{Preliminaries}

\subsection{Definitions}

By $\bE^d$ we denote the \emph{$d$-dimensional Euclidean space}
with \emph{origin} $O$. We write capital letters, e.g. $A_i$, for
\emph{points} or their \emph{position vectors}, and small letters,
like $a_i$, for their \emph{barycentric coordinates} in a suitably
defined system. Also we write $\|A_i - A_j \|$ for the usual
\emph{distance} between the points $A_i$ and $A_j$.

A (non-degenerate) \emph{$d$-simplex} $S = [A_1, \dots, A_{d+1}], \, d \ge 2$, is defined as
the convex hull of $d+1$ affinely independent points (or position vectors) $A_1, \dots, A_{d+1}$
in the Euclidean space $\bE^d$. The points $A_i$ are the \emph{vertices}, the line segments
$E_{ij} = A_i A_j$ (joining two different vertices $A_i, A_j$) the \emph{edges}, and all
$k$-simplices whose vertices are $k+1$ vertices of $S$ the \emph{$k$-faces} of $S$. The \emph{facets}
of $S$ are its $(d-1)$-faces, and the \emph{i-th facet} $F_i$ is the \emph{facet opposite to
the vertex} $A_i$.

The \emph{centroid} $\cG$ of $S$ is the average of its vertices,
the \emph{circumcenter} $\cC$ is the center of the unique sphere
containing all vertices, and the \emph{incenter} $\cI$ is the
center of the unique sphere that touches all facets of $S$. The
corresponding radii are called \emph{circumradius} $R$ and
\emph{inradius} $r$, respectively. For $d=2$, the altitudes of $S$
have a common point, the \emph{orthocenter} $\cH$ of $S$. For $d
\ge 3$ the altitudes of $S$ are not necessarily concurrent, but if
they are, i.e., if they have the orthocenter $\cH$ in common, we
say that $S$ is an \emph{orthocentric $d$-simplex}. A special
class of orthocentric $d$-simplices, also investigated here, is
that of \emph{rectangular $d$-simplices}; in this case the
orthocenter coincides with a vertex of $S$. The last center to be
defined here exists again for any $d$-simplex. Namely, for each
edge $E_{ij}$ of $S$ there is a unique hyperplane $H_{ij}$
perpendicular to $E_{ij}$ and containing the centroid $\cG_{ij}$
of the remaining $d-1$ vertices. These ${d+1 \choose 2}$
hyperplanes have a common point, the \emph{Monge point} $\cM$ of
$S$. This point is the reflection of $\cC$ in $\cG$ and coincides,
if $S$ is orthocentric, with the orthocenter. For the history of
the Monge point we refer to \cite{Za}, \S~21, and \cite{Cr}.

A $d$-simplex $S$ is said to be \emph{regular} or
\emph{equilateral} if all its edges have the same length. Note
that this is the highest degree of symmetry that $S$ can have, in
the sense that the group of isometries of a regular $d$-simplex is
the full symmetric group of permutations on the set of its
vertices, see \cite{Be}, Proposition 9.7.1. Several weaker degrees
of regularity of simplices are discussed in the present paper,
too. In particular, a $d$-simplex $S$ is said to be
\emph{equifacetal} if all its facets are congruent (or isometric),
and it is called \emph{equiareal} if all its facets have the same
$(d-1)$-volume, i.e., $(d-1)$-dimensional Lebesgue measure.
Furthermore, a $d$-simplex satisfying $\cG = \cC = \cI$ may be
referred to as $(\cG, \cC, \cI)$--\emph{equicentral}. (Note that
for $d=2$ the latter three degrees of regularity are equivalent to
equilaterality.)

\subsection{Some results on centers of general simplices}

Interesting results on equifacetal, equiareal and related
simplices are contained in the papers \cite{De 1}, \cite{De 2},
\cite{F-M}, \cite{McM}, \cite{We}, \cite{M-W}, \cite{Ed}, and
\cite{E-H-M}.
 For $d \ge 2$ we have
\begin{eqnarray}
S \, \mbox{ is regular } \, & \Longrightarrow & S \, \mbox{ is equifacetal } \, \Longrightarrow
S \, \mbox{ is } \, (\cG, \cC, \cI)\!-\!\mbox{equicentral}\,, \\
\cI = \cG \, & \Longleftrightarrow &  S \, \mbox { is equiareal},
\end{eqnarray}
see \cite{E-H-M}, Theorem 3.2. The other coincidences $\cG = \cC$
and $\cC = \cI$ turn out to have other geometric interpretations
that are worth recording. Thus, calling a $d$-simplex
\emph{equiradial} if all its facets have equal circumradii, and
\emph{of well-distributed edge-lengths} (or \emph{equivariant}) if
the sum of the squared edge-lengths is the same for all facets, it
follows from \cite{E-H-M}, Theorem 3.2, that for $d \ge 2$
\begin{eqnarray}
\cG = \cC & \Longleftrightarrow & S \, \mbox{ has well-distributed edge-lengths },\\
\cC = \cI & \Longleftrightarrow & \cC \, \mbox{ is interior and } S \mbox{ is
equiradial }.
\end{eqnarray}

\subsection{The Gram matrix of a simplex}

We continue with the representation of a tool that relates the
geometry of a simplex $S$ to the algebraic properties of a certain
matrix associated to $S$ (see \cite{H-J} and \cite{L-T}). Namely,
for a $d$-simplex $S = [A_1, \dots, A_{d+1}]$ in $\bE^d$ one
defines the \emph{Gram matrix} $G$ to be the symmetric, positive
semidefinite $(d+1) \times (d+1)$ matrix of rank $d$ whose
$(i,j)$-th entry is the inner product $A_i \cdot A_j$ (we mean the
ordinary inner product, say), cf. \cite{H-J}, p. 407. Given $G$,
one can calculate the distances $\|A_i - A_j\|$ for every $i,j$
using the formula
\[
\left(\|A_i - A_j\|\right)^2 = \left(A_i - A_j\right) \cdot
\left(A_i - A_j\right)\,.
\]
According to the last part of Proposition 9.7.1 in \cite{Be}, $G$
determines $S$ up to an isometry of $\bE^d$. Also one recovers $S$
from $G$ via the \emph{Cholesky factorization} $G = HH^t$, where
the rows of $H$ are the vectors $A_i$ coordinatized with respect
to some orthonormal basis of $\bE^d$. In fact, if $G$ is a
symmetric, semidefinite, real matrix of rank $r$, say, then there
exists a unique symmetric, positive semidefinite, real matrix of
rank $r$ with $H^2 = G$, cf. \cite{H-J}, Theorem 7.2.6, p. 405,
and the symmetry of $H$ implies $G = HH^t$.

\section{Basic properties of orthocentric simplices and a survey of known results}

We start this section with a short survey of known results about
orthocentric $d$-simplices. Since the case $d=3$ is not in our
focus here, we mention only some basic references referring to
orthocentric (and closely related) tetrahedra, namely \cite{Za},
\S~21, \cite{Si}, \S~30, \cite{AC 1}, \cite{AC 2}, \cite{CB},
\cite{The}, \cite{AC}, Chapters IV and IX, and the recent paper
\cite{H-W}. (Even J. L. Lagrange \cite{La} obtained results about
orthocentric tetrahedra over 200 years ago.) So the following
short survey refers to results on orthocentric $d$-simplices for
all dimensions $d \ge 3$.

A first basic property of orthocentric simplices is the fact that
they are closed under passing down to faces. (This sometimes
allows one to use induction on the dimension for establishing
certain properties in high dimensions.) More precisely, \emph{each
$k$-face of an orthocentric $d$-simplex is itself orthocentric},
$2 \le k < d$. Even more, in this passing-down procedure the feet
of all altitudes of any $(k+1)$-face $F$ are the orthocenters of
the $k$-faces of $F$. These observations were often rediscovered,
see \cite{Lo}, \cite{E 1}, \cite{K-T}, \cite{R-J}, \cite{K-H-M},
and \cite{PT}, \S~1.3, Problems 1.28, 1.29 and their solutions.

Another fundamental property of orthocentric simplices is the
perpendicularity of non-intersecting edges. It can also be
formulated as follows: \emph{Each edge of an orthocentric
$d$-simplex, $d \ge 3$, is perpendicular to the opposite
$(d-2)$-face, and any $d$-simplex with that property is
orthocentric}, cf. \cite{Lo}, \cite{K-T}, \cite{R-J},
\cite{K-H-M}, \cite{B-W}, and \cite{PT}, \S~1.3, for various
approaches. Also in \cite{PT}, \S~1.3, one can find the following
property of an orthocentric simplex $S = [A_1, \dots, A_{d+1}]$
regarding its circumcenter: $\overrightarrow{\cC A_1} + \dots +
\overrightarrow{\cC A_{d+1}} = (d-1) \overrightarrow{\cC \cH}$
~(cf. Problem 1.29 there). One of the oldest and most elegant
discoveries in the geometry of triangle centers is  Euler's proof
that the centroid $\cG$ of a triangle $ABC$ lies on the line
segment $\cC \cH$ and divides it in the ratio $1:2$, see
\cite{Co}, p. 17. For an orthocentric $d$-simplex $S$ we have the
analogous situation: \emph{The points} $\cC, \cG$ \emph{and $\cH$
are on a line $($the Euler line of $S$$)$, and $\cG$ divides the
segment $\cC \cH$ in the ratio} $(d-1):2$. This result and its
analogue for general $d$-simplices (where the Monge point $\cM$
replaces $\cH$) were also rediscovered several times, see
\cite{Me}, \cite{Lo}, \cite{E 1}, \cite{Ha}, \cite{K-T},
\cite{Mo}, and \cite{B-W} for different proofs and extensions.
Also in some other situations, theorems on orthocentric simplices
have analogues for general simplices if the missing point $\cH$ is
replaced by $\cM$, such as in the case of Feuerbach spheres
discussed in the sequel. Once more we mention that if a simplex is
orthocentric, then $\cM$ coincides with $\cH$, see \cite{Me} and
\cite{Mo}.

It is well known that the ${d+1 \choose k+1}$ centroids of all
$k$-faces, $k \in \{0, \dots, d-1\}$, of an orthocentric
$d$-simplex $S$ lie on a sphere, the \emph{Feuerbach $k$-sphere}
of $S$, see \cite{Sh} and \cite{K-T}. Also for general simplices,
H. Mehmke \cite{Me} investigated the case $k=d-1$: the center
$\cF$ of the respective Feuerbach sphere lies on the Euler line
through $\cG$ and $\cC$, and $\cG$ divides the segment $\cF \cC$
in the ratio $1:d$. The radius of that sphere is $\frac{1}{d}$
times the circumradius of $S$. If $S$ is orthocentric, then the
Feuerbach sphere also contains the feet of the altitudes, and it
divides the ``upper'' parts of the altitudes in the ratio
$1:(d-1)$, see also \cite{Ri}, \cite{Lo}, \cite{Man}, \cite{Mo},
\cite{Fr}, and \cite{Fi 1} for analogous results. In \cite{E 1},
\cite{E 2}, \cite{Ha}, \cite{K-T}, and \cite{Ge 2} the whole
sequence of all Feuerbach $k$-spheres of orthocentric
$d$-simplices is studied. All their corresponding centers $\cF_0,
\dots, \cF_{d-1}$ lie on the Euler line, with $|\cH \cF_k| : |\cH
\cG| = (d+1) : 2 (k+1)$, and their radii $r_k$ satisfy simple
relations depending only on $d$ and $k$. These papers contain more
related results (see also \cite{Iw}), e.g. observations referring
to Feuerbach spheres of so-called orthocentric point systems.
Considering the set of $d+2$ points $\cH, A_1, \dots, A_{d+1}$ of
an orthocentric $d$-simplex as a whole, one observes that each of
them is the orthocenter of the simplex formed by the $d+1$ others.
Thus it is natural to speak about \emph{orthocentric systems of
$d+2$ points}. Such point sets (and their analogues of larger
cardinality) were studied in \cite{Ri}, \cite{E 1}, \cite{Fi 1},
and \cite{B-W}. For example, E. Egerv\'ary \cite{E 1} proved that
\emph{a point set $\{P_0, P_1, \dots, P_{d+1}\} \subset \bE^d\,,
\, d \ge 2$, is an orthocentric system if and only if the mutual
distances $\|P_i - P_j \| \, (i,j = 0,1, \dots, d+1; i \not= j)$
can be expressed by $d+2$ symmetric parameters $\lambda_i$ in the
form}
\begin{equation}
\| P_i - P_j \|^2 = \lambda_i + \lambda_j\,,\,\,
\mbox{ with } \, \sum \limits^{d+1}_{i=0} \frac{1}{\lambda_i} = 0\,, \,\, \lambda_i + \lambda_j > 0 \,, \,\, i \not= j\,,
\end{equation}
see also \cite{PT}, \S~1.3, Problem 1.28. Based on this, Egerváry
showed that $d+1$ points in $\bE^d$, whose cartesian coordinates
are the elements of an \emph{orthogonal matrix}, form together
with the origin an orthocentric point system. If, conversely, the
``interior point'' of an orthocentric system of $d+2$ points is
identified with the origin, then an orthogonal matrix can be found
from which the coordinates of the remaining $d+1$ points can be
easily described. In \cite{B-W} it is shown that the $d+2$ Euler
lines of an orthocentric system $\{P_0, P_1, \dots, P_{d+1}\}
\subset \bE^d\,, \, d \ge 2$, have a common point, called its
\emph{orthic point}, and that the $d+2$ centroids as well as the
$d+2$ circumcenters of that set form again orthocentric systems,
both homothetic to $\{P_0, P_1, \dots, P_{d+1}\}$ with the orthic
point as homothety center.

M. Fiedler \cite{Fi 1} defines \emph{equilateral $d$-hyperbolas}
as those rational curves of degree $d$ which have all their $d$
asymptotic directions mutually orthogonal. Two such $d$-hyperbolas
are called \emph{independent} if both $d$-tuples of asymptotic
directions satisfy the following: In no $k$-dimensional (improper)
linear subspace $(k = 1, \dots, d-2)$, which is determined by $k$
directions from one of these $d$-tuples, more than $k$ asymptotic
directions of the other are contained. He proves that \emph{if
there are two independent equilateral $d$-hyperbolas both
containing a system of $d+2$ distinct points in $\bE^d$, then this
system is orthocentric, and every $d$-hyperbola containing this
system is equilateral}.

Another type of results refers to \emph{characterizations of
orthocentric simplices as extreme simplices} regarding certain
metrical problems, going back to J. L. Lagrange \cite{La} and W.
Borchardt \cite{Bo}, and connected with the symmetric parameters
in (5), see \cite{E 1} and \cite{Ge}. In the latter paper the
following (and further) results are shown: \emph{The maximum
$[$minimum$]$ volume of a $d$-simplex $S = [A_1, \dots, A_{d+1}]$
containing a point $Q$ and with prescribed distances $\|Q-A_i\|
\ge 0\,, \, i = 1,\dots, d+1$, is attained by an orthocentric
$d$-simplex. The maximum volume of a $d$-simplex $S$ with given
$(d-1)$-volumes of its facets is attained if $S$ is orthocentric}.
For getting these results, L. Gerber \cite{Ge} establishes some
purely geometric properties that orthocenters of orthocentric
simplices must have, e.g.: The point $\cH$ of an orthocentric
simplex lies closer to a facet than to the opposite vertex on all
except possibly the shortest altitude. Further geometric
properties of orthocentric simplices were derived in \cite{Fi 1},
\cite{Fi 2}, and \cite{Fi 3}, \S~7. Going back to the parameters
$\lambda_i$ in (5), Fiedler \cite{Fi 1} calls an orthocentric
simplex \emph{negatively orthocentric} if one of the $\lambda_i$'s
is negative, \emph{positively orthocentric} if all of them are
positive, and \emph{singularly orthocentric} if one is zero (the
first two cases together are called \emph{non-singularly
orthocentric}). He shows that \emph{a $d$-simplex, $d \ge 2$, with
dihedral interior angles $\phi_{ij}$ is non-singularly
orthocentric if and only if there exist real non-zero numbers
$c_i\,, \, i=1, \dots, d+1$, such that} $\cos \phi_{ij} = c_i c_j$
\emph{for all} $i,j$ \emph{with} $i \not= j$.

A \emph{reciprocal transformation} with respect to the simplex $S$
is such that for the homogeneous barycentric coordinates of the
image $X' = (x'_i)$ of a point $X = (x_i)$ the relation $x'_i =
c_i \times x_i$ holds, where the $c_i$'s are fixed non-zero
numbers. The \emph{harmonic polar} of $Y = (y_i)$ not contained in
any face of $S$ is the hyperplane with equation $\sum
\frac{x_i}{y_i} = 0$ in barycentric coordinates. In \cite{Fi 1} it
is proved that \emph{a $d$-simplex $S$ is orthocentric if and only
if there exists an interior point $P$ of $S$ such that for every
selfadjoint point $Q$ $($if different from $P$$)$ of the
reciprocal transformation, for which $P$ and the centroid of $S$
correspond, the line through $P$ and $Q$ is perpendicular to the
harmonic polar of $Q$ with respect to $S$. The point $P$ is then
the orthocenter of $S$}. The paper \cite{Fi 1} contains also a
number of theorems on natural generalizations of positively
orthocentric simplices. In the booklet \cite{Fi 3} the
\emph{polarity} of a point quadric (in a projective space) with
equation $\sum a_{ik} x_i x_k = 0$ and a dual quadric $\sum b_{ik}
\xi_i \xi_k$ are defined, as usual, by the condition $\sum a_{ik}
b_{ik} = 0$, cf. Def. 7.8 there. A point quadric in $\bE^d$ is
then called \emph{equilateral} (cf. Def. 7.9 in \cite{Fi 3}) if it
is a polar to the absolute dual quadric. In the case of
homogeneous barycentric coordinates the dual absolute quadric has
the equation $\sum q_{ik} \xi_i \xi_k = 0$, where the matrix
$(q_{ik})$ is the Gram matrix of the outward normals of the
$d$-simplex $S$ normalized so that the sum of the normals is the
zero vector. Fiedler \cite{Fi 3} proves that for a
\emph{non-singularly orthocentric $d$-simplex$, d \ge 2$, every
equilateral quadric containing all its vertices contains the
orthocenter as well, and that, conversely, every quadric
containing all vertices and the orthocenter is necessarily
equilateral}. And coming back to equilateral $d$-hyperbolas (see
above), he proves in \cite{Fi 1} the following: \emph{Suppose that
in a non-singularly orthocentric $d$-simplex $S\,, \, d \ge 2$,
the orthocenter is not contained in any hyperplane orthogonally
bisecting an edge. Then there exists exactly one equilateral
$d$-hyperbola containing all vertices, the orthocenter and the
centroid of $S$}.

More general classes of simplices which are still closely related
to orthocentric ones were studied by S. R. Mandan (see \cite{Ma 1}
and \cite{Ma 2}) and M. Fiedler (cf. \cite{Fi 1} and \cite{Fi 2}).
The simplices under consideration in the papers of Mandan have two
(or more) subsets of their set of altitudes, each subset having a
common point. And Fiedler \cite{Fi 1} obtains theorems on a family
of simplices having the class of positively orthocentric simplices
as subfamily. In \cite{Fi 2} he investigates the related class of
\emph{cyclic simplices}. Also the paper \cite{Ge 2} should be
mentioned here.

From the literature we also know theorems on special types of
orthocentric simplices, in particular on the subfamily of
\emph{regular simplices}. So we know that an orthocentric
$d$-simplex $S = [A_1, \dots, A_{d+1}]$ is regular if and only if
$\cC = \cG$ \cite{Mo}, if and only if $\cG = \cH$ \cite{Fr}, if
and only if $\cH$ coincides with the unique point minimizing $\sum
\limits^{d+1}_{i=1} \| X - A_i \|\,, \, X \in \bE^d$ (the
Fermat-Torricelli point of $S$, see \cite{PW}), and if $S$ is
equiareal \cite{Ge}. Furthermore, there exist also some results on
rectangular (or right) simplices as special orthocentric ones, see
again \cite{Ge}.

\section{The barycentric coordinates of the orthocenter and the Gram matrix of an orthocentric simplex}

\subsection{Barycentric coordinates and obtuseness}

In this section, we show that a non-rectangular orthocentric
simplex  can essentially be parametrized by the barycentric
coordinates of its orthocenter, and we give several useful
characterizations of non-rectangular orthocentric simplices.

We start with a simple but basic theorem.

\textbf{Theorem 3.1:} \emph{Let $S = [A_1, \dots, A_{d+1}]$ be a
$d$-simplex.}
\begin{enumerate}
\item[(a)] \emph{$S$ is orthocentric if and only if for every $k$
the quantity $(A_i - A_k) \cdot (A_j - A_k)$ does not depend on
$i$ and $j$ as long as $i,j$ and $k$ are pairwise distinct.}
\item[(b)] \emph{If $P$ is a point in the affine hull of $S$, then
$S$ is orthocentric with orthocenter $P$ if and only if the
quantity $(A_i - P) \cdot (A_j - P)$ does not depend on $i$ and
$j$ as long as $i \not= j$.}
\end{enumerate}
\emph{If $c_k$ denotes the quantity in} (a) \emph{and $c$ the quantity in} (b), \emph{then we have
\begin{eqnarray}
c = 0 & \Longleftrightarrow & S \mbox{ is rectangular at } A_k \mbox{ for some } k\,.\\
c_k = 0 & \Longleftrightarrow & S \mbox{ is rectangular at } A_k \,.
\end{eqnarray}
}

\textbf{Proof:} For (a), see \cite{PT}, Problem 1.28, pages 30,
217. To prove (b), we may assume, without loss of generality, that
$P$ is the origin $O$. If $S$ is orthocentric with orthocenter
$\cH$, and if $i,j$, and $k$ are pairwise distinct indices, then
$A_i$ is normal to the $i$-th facet, and therefore $A_i \cdot (A_j
- A_k) = 0$, and $A_i \cdot A_j = A_i \cdot A_k$, as desired.
Conversely, if $A_i \cdot A_j$ does not depend on $i$ and $j$ as
long as $i \not= j$, then $A_i \cdot (A_j - A_k) = 0$, for all
pairwise distinct $i,j$, and $k$, and therefore $A_i$ is normal to
every edge of the $i$-th facet, and hence to the $i$-th facet.
Therefore $O$ lies on every altitude and has to be the
orthocenter. Finally, if $c = 0$, and if $\cH$ is the orthocenter
of $S$, then the $d+1$ vectors $\cH - A_1, \dots, \cH - A_{d+1}$
in $\bE^d$ are normal to each other, and therefore one of them
must be the zero vector, i.e., $\cH = A_i$ for some $i$. The other
implications follow from the definitions. \hfill $\Box$

\textbf{Definition 3.2:} For an orthocentric $d$-simplex $S =
[A_1, \dots, A_{d+1}]$ with orthocenter $\cH$, we define
 $\sigma (S)$ of $S$ by
\begin{equation}
\sigma (S): = (\cH - A_i) \cdot (\cH - A_j)\,, \,\, i \not= j\,.
\end{equation}
By Theorem 3.1(b), $\sigma (S)$ is well-defined. Note that in view
of (6), $\sigma (S)$ is zero  if and only if  $S$ is rectangular.
We shall see later that the sign of $\sigma (S)$ is negative if
and only if all the angles between any two edges of $S$ are acute.
Because of this and for lack of a better term, we propose to call
$\sigma (S)$ \emph{the obtuseness of} $S$. This may conceal the
fact that if $S$ is enlarged (by a factor of $\lambda$, say), then
also $\sigma (S)$ increases (by a factor of $\lambda ^2$),
although the shape of $S$ remains unchanged.

\bigskip
The next technical theorem will be freely used. Among other
things, it expresses the edge-lengths of an orthocentric
non-rectangular $d$-simplex $S$ in terms of the obtuseness $\sigma
(S)$ of $S$  and the barycentric coordinates of its orthocenter
$\cH$, showing that these quantities are sufficient for
parametrizing such simplices. Clearly, this does not apply to
rectangular simplices, since these numbers do not carry any
information at all about the simplex other than its being
rectangular. This is the main reason why rectangular simplices are
temporarily excluded and their study is postponed to a later
section.

\textbf{Theorem 3.3:} \emph{Let $S = [A_1, \dots, A_{d+1}]$ be a
non-rectangular orthocentric $d$-simplex, and let $a_1, \dots,
a_{d+1}$ be the barycentric coordinates of its orthocenter $\cH$
with respect to $S$. Let $c = \sigma (S)$ be the obtuseness of $S$
defined in $(8)$. Then no $a_i$ is equal to $0$ or $1$ and, for
any real numbers $b_1, \dots, b_{d+1}$,}
\begin{equation}
\left\| {\small \sum} b_i (A_i - \cH) \right\|^2 = c
\left[\left(\sum b_i\right)^2 - \sum \frac{b_i^2}{a_i}\right]\,.
\end{equation}
\emph{In particular},
\begin{eqnarray}
\|A_i - \cH \|^2 & = & \frac{c(a_i-1)}{a_i}\\[0.3cm]
\|A_i - A_j \|^2 & = & - c \left(\frac{1}{a_i} + \frac{1}{a_j}\right)\,.
\end{eqnarray}
\emph{Also, if $B_i$ is the foot of the perpendicular from the
vertex $A_i$ to the $i$-th facet, and if $h_i = \|A_i - B_i\|$ is
the corresponding altitude, then}
\begin{eqnarray}
B_i - \cH & = & \frac{a_i}{a_i -1} (A_i - \cH)\\[0.3cm]
h^2_i & = & \frac{c}{a_i (a_i-1)}\,.
\end{eqnarray}

\textbf{Proof:} Without loss of generality, we may assume that the orthocenter $\cH$ of $S$
lies at the origin $O$. Taking the scalar product of $A_i$ with $\sum a_i A_i = O$,  we see that $a_i \|A_i\|^2 + (1 - a_i) c = 0$. If $a_i = 0$, then $c = 0$.
If $a_i = 1$, then $\|A_i\| = 0$. In both cases $S$ would be rectangular. Therefore no $a_i$ is
$0$ or $1$, and
\[
\|A_i\|^2 = \frac{c(a_i-1)}{a_i}\,,
\]
as claimed in (10). It follows that
\[
\begin{array}{lll}
\|\sum b_i A_i\|^2 & = & \sum b^2_i \|A_i\|^2 + 2 \sum \limits_{i<j} b_i b_j (A_i
\cdot A_j)\\[0.3cm]
& = & c \left[\sum b^2_i \left({\displaystyle
\frac{a_i-1}{a_i}}\right) + 2 \sum \limits_{i<j}
b_i b_j \right]\\[0.5cm]
& = & c \left[\sum b^2_i - \sum {\displaystyle \frac{b^2_i}{a_i}} + 2 \sum \limits_{i<j}
b_i b_j \right]\\[0.5cm]
& = & c \left[\left(\sum b_i\right)^2 - \sum {\displaystyle \frac{b^2_i}{a_i}}\right]\,,
\end{array}
\]
as claimed in (9).

Next, let $\cH' = B_{d+1}$ be the projection of $A_{d+1}$ on the
$(d+1)$-th facet $F_{d+1}$ of $S$. Since $\cH'$ lies on the
(well-defined) line joining the vertex $A_{d+1}$ and the origin
$O$, and also in the affine hull of $F_{d+1}$, it follows that there
exist $t, a'_1, \dots, a'_d$ such that
\[
\cH' = t A_{d+1} \,\, \mbox{ and } \,\, \cH' = a'_1 A_1 + \dots + a'_d A_d\,,
\]
with $a'_1 + \dots + a'_d = 1$. Set $a'_{d+1} = -t$. Then we have
\[
a'_1 A_1 + \dots + a'_{d+1} A_{d+1} = O \,\, \mbox{ and } \,\, a_1 A_1 + \dots + a_{d+1}
A_{d+1} = O\,.
\]
From the uniqueness of the dependence relation among $A_1, \dots, A_{d+1}$ it follows that the
$(d+1)$-tuples $(a'_1, \dots, a'_{d+1})$ and $(a_1, \dots, a_{d+1})$ are proportional. Since
$a'_1 + \dots + a'_d = 1 \not= 0$, it follows that $a_1 + \dots + a_d \not= 0$, and that
\[
\frac{a'_i}{a'_1 + \dots + a'_d} = \frac{a_i}{a_1 + \dots + a_d}\,.
\]
Therefore
\begin{equation}
\frac{a'_i}{1} = \frac{a_i}{1 - a_{d+1}}\,, \,\, \cH' = -a'_{d+1} A_{d+1} =
\frac{a_{d+1}}{a_{d+1}-1} A_{d+1}\,.
\end{equation}
This proves (12) for $i = d+1$, and hence for all $i$. Finally,
\[
h^2_i = \|A_i - B_i\|^2 = \| A_i - \frac{a_i}{a_i-1} A_i \|^2 = \frac{c}{(a_i -1)^2}
\frac{a_i-1}{a_i} = \frac{c}{a_i (a_i-1)}\,,
\]
as claimed in (13). ~\hfill $\Box$

\bigskip

\textbf{Theorem 3.4:} \emph{Let $S = [A_1, \dots, A_{d+1}]$ be an
orthocentric $d$-simplex, and let $a_1, \dots, a_{d+1}$ be the
barycentric coordinates of its orthocenter $\cH$ with respect to
$S$. Let $\sigma (S)$ be the obtuseness of $S$ defined in $(8)$.
Then $S$ is non-rectangular if and only if any of the following
conditions hold.}
\begin{enumerate}
\item[(a)] $\sigma (S) \not= 0$.
\item[(b)] \emph{None of the faces of $S$ is rectangular}.
\item[(c)] \emph{$\cH$ is not in the affine hull of any proper face of $S$.}
\item[(d)] \emph{There does not exist any nonempty subset $I$ of $\{1, \dots, d+1\}$ such
that $\Sigma (a_i: i \in I ) = 0$.}
\item[(e)] \emph{There does not exist any proper subset $I$ of $\{1, \dots, d+1\}$ such that
$\Sigma (a_i : i \in I )= 1$.}
\end{enumerate}

\textbf{Proof:}  Without loss of generality, we may assume that
the orthocenter $\cH$ of $S$ lies at the origin $O$. Note that
property (a) has already been mentioned in (6). To prove (b),
suppose that $S$ has a rectangular face $F$. For simplicity, we
may assume that $F= [A_1, \dots, A_{k+1}], 2 \le k \le d$, and
that $A_1$ is the orthocenter of $F$. Then $(A_1 - A_2) \cdot (A_1
- A_3) = 0$. Since $S$ is orthocentric, it follows from Theorem
3.1 (a) that ($A_1 - A_i) \cdot (A_1 - A_j) = 0$ for all $i \not=
1$ and $j \not= 1$. Therefore the edges of $S$ emanating from
$A_1$ are normal to each other, and hence $S$ is rectangular. To
prove (c), suppose that the orthocenter of $S$ is in the affine
hull of a proper face $F$ of $S$. Without loss of generality, we
assume that $\cH = O$. If $A_i$ is a vertex of $F$, and $A_k$ is
not a vertex of $F$, then the segment $A_k \cH$, being normal to $
F$, is normal to $A_i \cH$. Thus $\sigma (S) = 0$, and $S$ would
be rectangular. This proves (c). To prove (d), suppose that $a_1 +
\dots + a_k = 0$ for some $k \ge 1$. By (9), we have $\|a_1 A_1 +
\dots +a_k A_k \| = 0$, and therefore $a_1 A_1 + \dots + a_k A_k =
O$ and $\cH = a_{k+1} A_{k+1} + \dots + a_{d+1} A_{d+1}$. Thus
$\cH$ is in the affine hull of the proper face $[A_{k+1}, \dots,
A_{d+1}]$, contradicting (c). This proves (d) and its equivalent
(e). Thus non-rectangular orthocentric simplices satisfy all the
conditions above. The converse is trivial. \hfill $\Box$

\subsection{The Gram matrix of an orthocentric simplex and
characterizing the barycentric coordinates of its orthocenter}

We now characterize those tuples that can occur as the barycentric
coordinates of a non-rectanglular orthocentric simplex, and we see
how the signs of these coordinates bear on the ``acuteness'' of
its vertex angles. We start with a theorem that describes the Gram
matrix (see our Subsection 1.3) of a non-rectangular orthocentric
simplex, and the lemma that follows records the value of a
determinant that we shall need in several places. This lemma
appears as Problem 192 on page 35 (with a hint on page 154, and an
answer on page 187) of \cite{F-S}.

\textbf{Theorem 3.5:} \emph{If $S = [A_1, \dots, A_{d+1}]$ is a
non-rectangular orthocentric $d$-simplex, then its Gram matrix is
of the form $cG$, where every entry of $G$ that lies off the
diagonal is $1$. Conversely, if $G$ is a $(d+1) \times
(d+1)$-matrix of rank $d$ such that its non-zero eigenvalues are
real and have the same sign, and such that all off-diagonal
entries are equal, say, the entries $G_{ij}$ of $G$ are given by}
\[
G_{ij} = \left\{\begin{array}{lll}
c(1+x_i) & \mbox{ if } & i = j\,,\\
c & \mbox{ if } & i \not= j\,,
\end{array}\right.
\]
\emph{then $\pm G$ is the Gram matrix of a non-rectangular
orthocentric $d$-simplex $S$ whose orthocenter lies at the origin
and for which
\begin{equation}
\|A_i\|^2 = c(1+x_i)\,, \,\, \|A_i - A_j\|^2 = c(x_i + x_j) \,\,
\,\,  \mbox{ if } \,\, \,\, i \not= j\,.
\end{equation}
}
 \textbf{Proof:} This is a restatement of the fact that $[A_1,
\dots, A_{d+1}]$ is orthocentric with orthocenter $O$ if and only
if $A_i \cdot A_j$ is a constant $c$ independent of $i$ and $j$
for $i \not= j$, and that $c = 0$ if and only if $S$ is
rectangular. \hfill $\Box$

\textbf{Lemma 3.6:} \emph{Let $\bJ = \bJ (a_1, \dots, a_n; b_1,
\dots, b_n)$ be the $n \times n$ matrix whose $(ij)$-th entry
$\bJ_{ij}$ is defined by
\[
\bJ_{ij} = \left\{ \begin{array}{lll}
a_i + b_i & \mbox{ if } & i = j\,,\\
a_i & \mbox{ if } & i \not= j\,.
\end{array} \right.
\]
Then}
\[
\det (\bJ (a_1, \dots, a_n; b_1, \dots, b_n)) = (b_1 b_2 \dots b_n)
\left( 1 + \frac{a_1}{b_1} + \frac{a_2}{b_2} + \dots + \frac{a_n}{b_n}\right)\,,
\]
\emph{where the right hand side is appropriately interpreted if $b_i = 0$ for some $i$. In particular}
\[
\det (\bJ (a, \dots, a; b, \dots, b)) = b^{n-1} (b+na)\,.
\]

{\bf Proof.} We proceed by induction, assuming that the statement is true for $n=k$.
Since det is continuous, we may confine ourselves to the domain where no $b_i$ is $0$.
By dividing the $i$--th row by $b_i$ for every $i$, we may also assume that $b_i=1$ for all $i$.
Then
\begin{eqnarray*}
&\det&\!\!\!\!(\bJ (a_1, \cdots, a_{k+1},; 1, \cdots, 1))\\[0.3cm]
& = & \left| \begin{array}{ccccc}
1+a_1 & a_1 & a_1 & \cdots & a_1\\
a_2 & 1+a_2 & a_2 & \cdots & a_2\\
\vdots & \vdots & \vdots & \vdots & \vdots\\
a_{k+1} & a_{k+1} & a_{k+1} & \cdots & 1+a_{k+1} \end{array}\right| \\[0.3cm]
& = & \left| \begin{array}{ccccc}
a_1 & a_1 & a_1 & \cdots & a_1\\
a_2 & 1+a_2 & a_2 & \cdots & a_2\\
\vdots & \vdots & \vdots & \vdots & \vdots\\
a_{k+1} & a_{k+1} & a_{k+1} & \cdots & 1+a_{k+1} \end{array}\right| +
\left| \begin{array}{ccccc}
1 & a_1 & a_1 & \cdots & a_1\\
0 & 1+a_2 & a_2 & \cdots & a_2\\
\vdots & \vdots & \vdots & \vdots & \vdots\\
0 & a_{k+1} & a_{k+1} & \cdots & 1+a_{k+1} \end{array}\right|\\[0.3cm]
& = & \left| \begin{array}{ccccc}
a_1 & 0 & 0 & \cdots & 0\\
a_2 & 1 & 0 & \cdots & 0\\
\vdots & \vdots & \vdots & \vdots & \vdots \\
a_{k+1} & 0 & 0 & \cdots & 1 \end{array}\right| + \det (\bJ(a_2, \cdots, a_{k+1}; 1, \cdots, 1))\\[0.3cm]
& = & a_1 + (1+a_2 + \cdots + a_{k+1}) = 1 + a_1 + a_2 + \cdots + a_{k+1}, \,\, \mbox { as desired }.
\qquad \qquad \Box
\end{eqnarray*}

\bigskip
The next theorem describes the somewhat surprising restrictions
that the barycentric coordinates $a_1,\dots,a_{d+1}$ of the
orthocenter of a non-rectangular orthocentric simplex $S$ must
obey. We find it more convenient to include in it the relations
among the signs of $a_1,\dots,a_{d+1}$, the sign of the obtuseness
$\sigma (S)$, and the acuteness of the angles between the edges.
Note in particular that the sign of the obtuseness $\sigma (S)$ of
$S$ can be read off the barycentric coordinates of the orthocenter
of $S$.

\textbf{Definition 3.7:} A polyhedral angle with vertex at $O$ and
with ``arms'' $OV_1, \dots, OV_n$ is called \emph{strongly acute}
(respectively, \emph{strongly obtuse, right}) if and only if
$\angle V_i O V_j$ is acute (respectively, obtuse, right) for all
$i \not= j$.

\textbf{Theorem 3.8:} \emph{The real numbers $a_1, \dots, a_{d+1}$
with $a_1 + \cdots + a_{d+1} = 1$ occur as the barycentric
coordinates of the orthocenter of a non-rectangular orthocentric
$d$-simplex if and only if all of them are positive, or exactly
one of them is positive and the others are negative. In the first
case $\sigma (S) < 0$, and all the vertex angles of $S$ are
strongly acute. In the second case $\sigma (S) > 0$, and one
vertex angle of $S$ is strongly obtuse while the others are
strongly acute}.

\textbf{Proof:} Let $S = [A_1, \dots, A_{d+1}]$ be a
non-rectangular orthocentric $d$-simplex, and $a_1, \dots,
a_{d+1}$ be the barycentric coordinates of the orthocenter $\cH$
of $S$. Assume that $\cH$ is at the origin $O$. We use the facts
that if $i \not= j$, then
\[
A_i \cdot A_j = c\,, \,\,\, \|A_i\|^2 = \frac{c(a_i-1)}{a_i}\,, \,\,\, \|A_i - A_j\|^2 = - c
\left(\frac{1}{a_i} + \frac{1}{a_j}\right)\,.
\]
If $c < 0$, then $a_i (a_i-1)<0$ for all $i$, and therefore $0 < a_i < 1$ for all $i$. In this case,
we have for pairwise distinct $i,j$, and $k$
\[
(A_i - A_j) \cdot (A_i - A_k) = \frac{c(a_i-1)}{a_i} - c - c + c = \frac{-c}{a_i} > 0\,.
\]
If $c > 0$, and if $a_1$ and $a_2$ are positive, then $\|A_1 - A_2\|^2$ would be negative. Therefore at
most one of the $a_i$'s is positive. Since the sum is $1$, it follows that exactly one $a_i$ is positive.
Here again, if $i,j$, and $k$ are pairwise distinct, and if $t_i = (A_i - A_j) \cdot (A_i - A_k)$, then $t_i
< 0$ if $a_i > 0$, and $t_i > 0$ if $a_i < 0$.

Conversely, suppose that $a_1, \dots, a_{d+1}$ are non-zero real
numbers whose sum is $1$. Let $x_i = -1/a_i$, and let $G$ be the
$(d+1) \times (d+1)$-matrix whose $(ij)$-th entry $G_{ij}$ is
given by
\[
G_{ij} = \left\{\begin{array}{lll}
1 + x_i & \mbox{ if } & i = j\,,\\
1 & \mbox{ if } & i \not= j\,.
\end{array} \right.
\]
By Lemma 3.6, the characteristic polynomial of $G$ is given by
\[
\begin{array}{lll}
F(\lambda) & = & (\lambda-x_1) \cdots (\lambda - x_{d+1}) \left(1 + {\displaystyle \frac{-1}{\lambda - x_1}}
+ \dots + {\displaystyle \frac{-1}{\lambda - x_{d+1}}}\right)\\
& = & f(\lambda) - f' (\lambda) \,,
\end{array}
\]
where
\[
f(\lambda) = (\lambda - x_1) \cdots (\lambda - x_{d+1})\,.
\]
Let $g(\lambda) = - e^{-\lambda} f (\lambda)$. Then $g' (\lambda)
= e^{-\lambda} (f(\lambda)-f'(\lambda)) = e^{-\lambda}
F(\lambda)$. The assumption $a_1 + \cdots + a_{d+1} = 1$ is
equivalent to saying that $g'(0) = 0$.

If all the $a_i$'s are
positive, then all the $x_i$'s are negative, and therefore $g$ has
$d+1$ negative zeros. Hence $g'$ has $d$ negative zeros. Thus $d$
of the zeros of $g'$ are negative and the remaining one is $0$.
Therefore $F$ has the same property, and hence $- G$ is the Gram
matrix of some $d$-simplex $S$. It is easy to see that such a
simplex $S$ is orthocentric with orthocenter $O$ and with
barycentric coordinates of $O$ as desired.

If exactly one of the $a_i$'s is positive, say
\[
a_1 \le \dots \le a_d < 0 < a_{d+1}\,,
\]
then
\[
x_1 \ge \dots \ge x_d > 0 > x_{d+1}\,,
\]
and $g$ has $d$ positive zeros $x_1, \dots, x_d$. Therefore $g'$
has $d-1$ positive zeros in $(0,x_1)$. Also, $g'$ has a zero in
$(x_1, \infty)$ since $g(x_1) = g(\infty) = 0$. Thus $d$ of the
zeros of $g'$ are positive, and the remaining one is $0$.
Therefore $F$ has the same property, and hence $G$ is the Gram
matrix of some $d$-simplex $S$. It is easy to see that such an $S$
is orthocentric with orthocenter $O$ and with barycentric
coordinates of $O$ as desired. \hfill $\Box$

\textbf{Corollary 3.9:} \emph{If a $d$-simplex $S = [A_1, \dots,
A_{d+1}]$ is orthocentric, then at least $d$ of its vertex
polyhedral angles are strongly acute, while the remaining one is
either right, strongly acute, or strongly obtuse}.

This gives us a seemingly clear-cut idea of which orthocentric
simplices ought to be called acute and which ought to be called
obtuse. On the other hand, we do not get the criterion that the
orthocenter is interior if and only if the circumcenter is.

\textbf{Theorem 3.10:} \emph{Let $S = [A_1, \dots, A_{d+1}]$ be an
orthocentric $d$-simplex. Then the circumcenter $\cC$ of $S$ is
interior if and only if $S$ is non-rectangular and if the
barycentric coordinates $a_1, \dots, a_{d+1}$ of the orthocenter
$\cH$ of $S$ are such that $0 < a_i < 1/(d-1)$ for all $i$.
Consequently, if the circumcenter is interior, then so is the
orthocenter, but not conversely}.

\textbf{Proof:} We will see in Section 5 that the circumcenter of
a rectangular $d$-simplex $S$ lies on the hypotenuse facet if
$d=2$, and lies outside $S$ if $d > 2$. So let $S = [A_1, \dots,
A_{d+1}]$ be non-rectangular orthocentric, and let $a_1, \dots,
a_{d+1}$ be the barycentric coordinates of the orthocenter $\cH$
of $S$. Assume that $\cH = O$. Then it follows from the Euler line
theorem (see, e.g., \cite{E 1} and \cite{Ha}) that the
circumcenter $\cC$ of $S$ is given by
\[
\cC = \frac{d+1}{2}\,\cG = \frac{1}{2} (A_1 + \dots + A_{d+1})\,,
\]
and therefore
\[
\cC = \sum \left( \frac{1}{2} + \frac{(1-d)a_i}{2}\right)A_i\,.
\]
Thus the $i$-th barycentric coordinate of $\cC$ is
$(1+(1-d)a_i)/2$, and $\cC$ is interior if and only if $a_i <
1/(d-1)$. It remains to show that the $a_i$'s must be all
positive. If not, then exactly one of them is positive (and less
than $1/(d-1)$) and the others are negative, contradicting the
fact that $a_1 + \dots + a_{d+1} = 1$. \hfill $\Box$

\section{Coincidence of centers of non-rectangular orthocentric simplices}

\subsection{Center coincidences except circumcenter = incenter}

We start with a simple application of the fact that orthocentric
simplices are closed under passing down to faces. Note that the
part of Theorem 4.1 connected with the assumption $\cH = \cG$ is
known, see \cite{Fr}. Note also that a stronger result will be
established later, in Theorem 4.3.

\textbf{Theorem 4.1:} \emph{Let $S$ be an orthocentric
$d$-simplex. If the orthocenter of $S$ coincides with the
circumcenter or with the centroid of $S$, then $S$ is regular}.

\textbf{Proof:} Being trivially true for a triangle, the statement immediately follows for
all $d$ by induction, using the facts that the circumcenter and the orthocenter of a
facet of $S$ are the orthogonal projections of the respective centers of $S$ on that
facet, and that the centroid and the orthocenter of a facet of $S$ are the intersections
of the respective cevians of $S$ with the facet. \hfill $\Box$

\textbf{Remark 4.2:} The proof of Theorem 4.1 above does not work
if the circumcenter or the centroid is replaced by the incenter,
since the incenter of a facet $F$ of a simplex $S$ is neither the
projection of the incenter of $S$ on the facet, nor the
intersection of the respective cevian through the incenter of $S$
with that facet. However, we shall see in Theorem 4.3 that the
statement itself remains valid also in that case. \hfill $\Box$

The following theorem will be strengthened later to include the
case $\cC = \cI$ and to include rectangular simplices. Note that
the part pertaining to $\cH = \cG$ is treated in \cite{Fr}, and
that the part pertaining to $\cG = \cC$ is treated in  \cite{Mo}.
However, we include them in order to give a unified approach.

\textbf{Theorem 4.3:} \emph{Let $S$ be a non-rectangular
orthocentric $d$-simplex. If $\cH = \cG, \cH = \cC, \cH = \cI, \cG
= \cI$, or $\cG = \cC$, then $S$ is regular.}

\textbf{Proof:} Let $a_1, \dots a_{d+1}$ be the barycentric
coordinates of the orthocenter $\cH$ of $S$ with respect to $S$
and assume, without loss of generality, that $\cH$ lies at the
origin $O$. We freely use the facts that $A_i \cdot A_j = c$ for
all $i \not= j$ and $\|A_i\|^2 = c (a_i - 1) / a_i$ for all $i$.
The cases $\cH = \cG$ and $\cH = \cC$ are dealt with in Theorem
4.1. Alternatively, if $\cH = \cG$, then the $a_i$'s are all
equal, and therefore the edge-lengths are equal, by (11). If $\cH
= \cC$, then $\|A_i\|$ does not depend on $i$, and therefore
$\|A_i - A_j\|$ does not depend on $i$ and $j$, as long as $i
\not= j$, and $S$ is regular. If $\cH = \cI$, then $a_i$ is
proportional to the $(d-1)$-volume of the $i$-th facet which in
turn is inversely proportional to the altitude $h_i$ from the
$i$-th vertex. Thus $a_i h_i$ is independent of $i$. From (13) we
have
\[
a^2_i h^2_i = \frac{a^2_i c}{a_i (a_i-1)} = \frac{ca_i}{a_i -1}\,.
\]
Since $x/(x-1)$ is $1-1$, it follows that the $a_i$'s are equal, and $S$ is again regular.

If $\cG = \cI$, then $S$ is equiareal by (\cite{E-H-M}, Theorem
3.2 (iii)), and therefore the altitudes are equal. Hence $a_i
(a_i-1) = a_j (a_j - 1)$ for all $i$ and $j$. Therefore $(a_i -
a_j) (a_i + a_j - 1) = 0$. By Theorem 3.4 (e), $a_i + a_j$ cannot
be equal to $1$. Therefore the $a_i$'s are all equal, and $S$ is
regular. (This case is also treated in \cite{PW}.)

If $\cG = \cC$, then $S$ has well-distributed edge-lengths (by \cite{E-H-M}, Theorem 3.2 (i)). Using
(11), one can easily see that the sum of the squares of the edge-lengths of the $i$-th facet is
given by
\[
c(d-1) \left[\left(\frac{1}{a_1} + \dots + \frac{1}{a_{d+1}}\right)-\frac{1}{a_i}\right]\,.
\]
Therefore, having well-distributed edge-lengths is equivalent to saying that the $a_i$'s are equal,
and that $S$ is regular. Alternatively, if $\cC = \cG$, then the circumcenter is given by
\[
\cC = \sum \limits^{d+1}_{i=1} \frac{1}{d+1} A_i\,,
\]
and $\|\cC-A_i\|$ does not depend on $i$. Then
\[
\begin{array}{lll}
(d+1) \| \cC - A_i \|^2
& = & \| (d+1) \cC - (d+1) A_i \|^2\\
& = & (d+1)^2 \|\cC\|^2 + (d+1)^2 \|A_i\|^2 - 2 (d+1) \sum \limits^{d+1}_{j=1} A_i \cdot A_j\\[0.3cm]
& = & (d+1)^2 \|\cC\|^2 + (d+1)^2 \|A_i\|^2 - 2 (d+1) (\|A_i \|^2 + dc)\\[0.3cm]
& = & (d+1)^2 \|\cC \|^2 + (d^2 - 1) \| A_i \|^2 - 2 (d+1) dc\,.
\end{array}
\]
Thus $\|A_i\|$ is independent of $i$, and $O$ is the circumcenter
of $S$. Therefore $\cH = \cC$, and $S$ is regular.  This completes
the proof. \hfill $\Box$

\subsection{Equiradial orthocentric simplices and kites}

In this subsection we begin a general study of equiradial
orthocentric simplices. In view of the fact that
\[
\cC = \cI \, \mbox{ if and only if } \, \cC \, \mbox{ is interior
and the simplex is equiradial } \,,
\]
it is natural to study the broader class of equiradial simplices
and then to single out those for which $\cC$ is interior.

It is convenient to give a name to those $d$-simplices in which
$d$ vertices form a regular $(d-1)$-simplex $T$, called the base, and are
at equal distances from the remaining vertex, called the apex. We propose
to call such a $d$-simplex a $d$-\emph{dimensional kite}, or simply a
$d$-\emph{kite}, and to denote it by $K_d [s,t]$, where $s$ is the
side-length of each edge of $T$, and where $t$ is the length of each
remaining edge. The subscript $d$ may be omitted if no confusion is
caused. The quotient $t/s$ carries all the information about the shape of
the kite, and will be called the \emph{eccentricity} of the kite. Note
that a kite is automatically orthocentric, since the altitude from its
apex meets the base in its orthocenter. Note also that all the facets of a
kite are themselves kites.

\textbf{Lemma 4.4:} \emph{
For $d\le 3$ an equiradial orthocentic $d$-simplex is regular. For
$d > 3$ there exists a unique siimilarity class of equiradial
$d$-simplices $S$ such that $d$ of its vertices
form a regular $(d-1)$-simplex $T$ (of edge-length $s$, say) and are
at the same distance
($t$, say) from the remaining vertex. Its eccentricity is given by
\begin{equation}
\epsilon=t/s = \sqrt{(d-2)/d}\,.
\end{equation}
}

\textbf{Proof:} For $d=2$, the result is trivial, so suppose $d\ge
3$. Suppose there is a non-regular equiradial $d$-kite $K$ with base
$T$. Then it has to  arise in the following way. Inscribe the regular
$(d-1)$-simplex $T = [A_1, \dots, A_d]$ in a $(d-2)$-hypersphere centered
at the origin in $\bE^{d-1}$, and having radius $\rho$, the circumradius
of $T$. Each facet of $T$ is the base of exactly two $(d-1)$-kites of
circumradius $\rho$: $T$ itself, and one other, which we now describe.
Let $A^*_i$ be the point diametrically opposite to $A_i$ for $1 \le i \le
d$. Thus $A^*_i = -A_i$. Let $s = \|A_i - A_j\|$ and $t = \|A^*_i -
A_j\|$, where $i \not= j$. By taking the inner product of each side of the
equation $A_1 + \cdots + A_d = O$ with itself, we obtain $\rho^2 d +
d(d-1) (A_i \cdot A_j) = 0$, and therefore $A_i \cdot A_j =
-\rho^2/(d-1)$, and
\begin{eqnarray}
s^2 & = & \|A_i - A_j\|^2 = 2\rho^2
\left(1 + \frac{1}{d-1}\right) = 2\rho^2
\left(\frac{d}{d-1}\right)\,,\\[0.3cm] t^2 & = & \|A^*_i - A_j\|^2 =
2\rho^2 \left(1 - \frac{1}{d-1}\right) = 2\rho^2
\left(\frac{d-2}{d-1}\right)\,.
\end{eqnarray}
It is clear that
\[
t = \rho \, \Longleftrightarrow \, d = 3 \,\, \mbox{ and that } \,\, t > \rho
\, \Longleftrightarrow \, d > 3\,.
\]
But if $t=\rho$, then $K$
degenerates, with the apex lying at the circumcenter of $T$. On the other
hand, if $t>\rho$, then an actual equiradial $d$-kite of positive height
can be formed from $T$ and the kites over the facets of $T$ and the above
calculations of $s^{2}$ and $t^{2}$ yield  the indicated formula for the
eccentricity. \hfill $\Box$

In view of the above
calculations, if $K$ is a $d$-dimensional kite whose regular base
has edge-length $s$ and circumradius $\rho$, then $\rho^2/s^2 =
(d-1)/(2d)$, and therefore the eccentricity of a $d$-kite can take
any value larger than $\sqrt{(d-1)/(2d)}$.

\textbf{Lemma 4.5:} \emph{A kite $K [s,t]$ in which the
circumcenter coincides with the incenter must be regular $($i.e.,
$ s = t$$)$}.

\textbf{Proof:} If $h$ is the altitude of a $d$-kite $K = K [s,t]$
to its regular base $T$, and if $\rho$ is the circumradius of $T$,
then, by Lemma 4.4,
\[
\begin{array}{lll}
&& \mbox{ the circumcenter of } K \mbox{ is interior }\\[0.3cm]
& \Longleftrightarrow & h^2 > \rho^2 \Longleftrightarrow t^2 - \rho^2 > \rho^2
\Longleftrightarrow t^2 > 2\rho^2\\[0.3cm]
& \Longleftrightarrow & t^2 > {\displaystyle \frac{d-1}{d}} s^2 \mbox{ (by (17))}\\[0.3cm]
& \Longleftrightarrow & {\displaystyle \frac{t^2}{s^2} > \frac{d-1}{d}}\,.
\end{array}
\]
By (18), the eccentricity of the non-regular equiradial $d$-kite
is given by $\sqrt{(d-2)/d}$, which is less than $\sqrt{(d-1)/d}$.
Therefore such kites do not have interior circumcenters, and their
circumcenters and incenters cannot coincide. Thus kites in which
$\cC = \cI$ must be regular. \hfill $\Box$

We record some basic formulas for quantities associated with a
kite. These formulas have fairly simple proofs essentially based
upon the Pythagorean Theorem, and they are also based on some
related calculations for the regular $d$-simplex which we give
first.

\textbf{Proposition 4.6:} \emph{Let $R=R_d=R_{d,s}$, $r=r_{d,s}$,
$h=h_{d,s}$, and $V=V_{d,s}$ denote, respectively, the
circumradius, the inradius, the altitude, and the volume of a
regular $d$-simplex of edge-length $s$. Then}
\begin{eqnarray}
R^2 &=& \frac{s^2 d}{2(d+1)} \label{R}\\
h &=& \sqrt{\frac{1}{2}} \sqrt{\frac{d+1}{d}} s \label{h}\\
V &=& \frac{1}{d !} \sqrt{\frac{d+1}{2^d}}      s^d\label{V} \\
r &=& \sqrt{\frac{1}{2d(d+1)}}~ s. \label{r}
\end{eqnarray}

\textbf{Proof:} To verify these formulas, let $S=[A_1,\cdots,A_{d+1}]$ be our
regular simplex, and assume that the center of $S$ is the origin
$O$. Then $A_1 + \cdots + A_{d+1} = O.$ Taking the scalar product
with $A_1$,  we see that $R^2 + d (A_1 \cdot A_2) = 0$, and
therefore $A_1 \cdot A_2 = -R^2/d.$ Also, $s^2 = (A_1 - A_2)^2 =
2R^2 + 2R^2/d.$ Therefore
\begin{eqnarray*} \label{sR-Rs}
    s^2 = \frac{2R^2(d+1)}{d}~,~~R^2 = \frac{s^2d}{2(d+1)}.
\end{eqnarray*}
This proves (\ref{R}).

By Pythagoras' Theorem, we have
\begin{eqnarray*}
    s^2 &=& h^2 + R_{d-1}^2 ~=~  h^2 + \frac{s^2(d-1)}{2d} \\
    h^2 &=& \frac{s^2(d+1)}{2d}.
    \end{eqnarray*}
This proves (\ref{h}).

Using (\ref{h}), we have
\begin{eqnarray*}
    V_d &=& \frac{1}{d}  h_d V_{d-1} ~=~
     \frac{1}{d}  \sqrt{\frac{1}{2}} \sqrt{\frac{d+1}{d}} s V_{d-1}~=~
     \frac{1}{d(d-1)}  \left(\sqrt{\frac{1}{2}}\right)^2 \sqrt{\frac{d+1}{d-1}} s^2 V_{d-2}\\ &=& ...... \\
     &=& \frac{1}{d !}  \left(\sqrt{\frac{1}{2}}\right)^{d-1} \sqrt{\frac{d+1}{2}} s^{d-1} V_{1} ~=~
      \frac{1}{d !}  \left(\sqrt{\frac{1}{2}}\right)^{d} \sqrt{d+1} s^{d}.
     \end{eqnarray*}
This proves (\ref{V}).

Finally, to calculate $r$, we use the fact that $V_d = (d+1) (r~
V_{d-1}/3)$ to obtain
\begin{eqnarray*}
r & = &   \frac{d}{d+1} \frac{V_d}{V_{d-1}} ~=~ \frac{d}{d+1}
\frac{1}{d}  \sqrt{\frac{d+1}{d}} \sqrt{\frac{1}{2}} ~ s~=~
\sqrt{\frac{1}{2d(d+1)}}~ s.
\end{eqnarray*}
~~~\hfill $\Box$

Now we are ready to prove the announced

\bigskip
\textbf{Theorem 4.7:} \emph{Let  $K$ be a $d$-kite whose regular
$(d-1)$-base $S$ has side-length $s$ and whose remaining vertex
$P$ is at distance $t$ from the vertices of $S$. Let $R=R_{d,s,t}$, $r=r_{d,s,t}$,
$h=h_{d,s,t}$, and $V=V_{d,s,t}$ denote, respectively, the
circumradius, the inradius, the altitude, and the volume of $K$.
Then}
\begin{eqnarray}
R &=& \sqrt{\frac{d}
{2\left(2\left(\frac{t}{s}\right)^2d - (d-1)\right)}}~  \left(\frac{t}{s}\right)^2 s \label{R-K}\\
h &=&  \sqrt{\frac{2\left(\frac{t}{s}\right)^2 d-(d-1)}{2d}} ~s\label{h-K}\\
V &=& \frac{1}{d !} \sqrt{\frac{2\left(\frac{t}{s}\right)^2 d-(d-1)}{2^d}}~s^{d} \label{V-K} \\
r &=&
 \frac{\sqrt{2\left(\frac{t}{s}\right)^2 d-(d-1)}}{\sqrt{2}
\left(\sqrt{d}+d\sqrt{2\left(\frac{t}{s}\right)^2 (d-1)-(d-2)} \right)}~s.
\label{r-K}
\end{eqnarray}

\textbf{Proof:} Let $u$ be the circumradius and $v$ the volume of
$S$. Then
\begin{eqnarray*}
    h^2 &=& t^2 - u^2 ~=~ t^2 - \frac{s^2(d-1)}{2d} ~=~
        \frac{2t^2d-s^2(d-1)}{2d}.
    \end{eqnarray*}
This proves (\ref{h-K}). On the other hand, $h = R \pm
\sqrt{R^2-u^2}$. Therefore,
\begin{eqnarray*}
t^2 - u^2  &=& \left( R \pm \sqrt{R^2-u^2} \right)^2
~=~ 2R^2 - u^2 \pm 2 R \sqrt{R^2-u^2}.\\
\left(t^2 - 2R^2\right)^2 &=&  4 R^2 (R^2-u^2)\\
t^4 - 4 R^2 t^2 &=& -4R^2 u^2   ~=~ -4 R^2 \frac{s^2(d-1)}{2d}\\
2  t^4 d&=&  4R^2 \left(2t^2d - s^2(d-1)  \right)\\
R^2 &=& \frac{t^4 d}{2\left(2t^2d - s^2(d-1)\right)}.
\end{eqnarray*}
This proves (\ref{R-K}). Also,
\begin{eqnarray*}
    V &=& \frac{1}{d} h v
~=~     \frac{1}{d} \sqrt{\frac{2t^2d-s^2(d-1)}{2d}}
\frac{1}{(d-1) !} \sqrt{\frac{d}{2^{d-1}}}~s^{d-1} \\
&=& \frac{1}{d !}    \sqrt{\frac{2t^2d-s^2(d-1)}{2^d}} ~
s^{d-1}.
    \end{eqnarray*}
This proves (\ref{V-K}). It remains to prove (\ref{r-K}). We use
\begin{eqnarray*}
     V &=&    V_{d,s,t} ~=~
     \frac{r}{d} v + d \left( \frac{r}{d} V_{d-1,s,t} \right) \\
     &=&
\frac{r}{d} \frac{1}{(d-1) !} \sqrt{\frac{d}{2^{d-1}}}~s^{d-1} + d
\left( \frac{r}{d} \frac{1}{(d-1) !}
\sqrt{\frac{2t^2(d-1)-s^2(d-2)}{2^{d-1}}}~s^{d-2}
\right) \\
     &=&
\frac{r}{d ! \sqrt{2^{d-1}}} s^{d-2}\left( \sqrt{d} s + d
\sqrt{2t^2(d-1)-s^2(d-2)} \right).
\end{eqnarray*}
Therefore
\begin{eqnarray*}
\frac{1}{d !} \sqrt{\frac{2t^2d-s^2(d-1)}{2^d}}~s^{d-1} &=&
\frac{r}{d ! \sqrt{2^{d-1}}} s^{d-2}\left( \sqrt{d} s + d
\sqrt{2t^2(d-1)-s^2(d-2)} \right),
\end{eqnarray*}
and hence
\begin{eqnarray*}
\sqrt{\frac{2t^2d-s^2(d-1)}{2}}~s &=& r \left( \sqrt{d} s + d
\sqrt{2t^2(d-1)-s^2(d-2)} \right).
\end{eqnarray*}
Therefore
\begin{eqnarray*}
r &=& s~ \frac{\sqrt{2t^2d-s^2(d-1)}}{\sqrt{2} \left( \sqrt{d} s +
d  \sqrt{2t^2(d-1)-s^2(d-2)} \right)},
\end{eqnarray*}
as desired.  \hfill $\Box$

\textbf{Addendum to Theorem 4.7:} The formulas derived above lead
to the following statements about the eccentricity $\epsilon =\frac{t}{s}$,
interior nature of the circumcenter, equiradiality, and the interior nature of the
orthocenter for a general kite:
\begin{quote}
$\epsilon^2$ can take any value in the interval $\left({\displaystyle \frac{d-1}{2d}}, \infty\right)$.\\
The circumcenter of $K$ is interior if and only if $\epsilon^2 > {\displaystyle \frac{d-1}{d}}$.\\
$K$ is equiradial if and only if $\epsilon^2 = {\displaystyle \frac{d-2}{d}}$ or
$\epsilon = 1$ (i.e., $t=s$).\\
The orthocenter of $K$ is interior if and only if $\epsilon^2
> {\displaystyle \frac{1}{2}}\,.$
\end{quote}

\subsection{Tools}

The following theorems are needed in the proof of Theorems 4.11 and
4.13. Theorem 4.8 describes how the barycentric coordinates $a_1,
\dots, a_{d+1}$ of the orthocenter of a non-rectangular
orthocentric $d$-simplex $S$ are related to those of the
orthocenter of a face $F$ of $S$, and also how $\sigma (F)$ and
$\sigma (S)$ are related. Theorem 4.9 expresses the circumcenter
of such an $S$ in terms of its vertices, and the circumradii of
$S$ and of its faces in terms of $a_1, \dots, a_{d+1}$ and $\sigma
(S)$. Lemma 4.10 deals with the Gram matrix of a special type of a
non-rectangular orthocentric $d$-simplex $S$. This is then used in
the proof of Theorem 4.11a,b which gives a characterization of
non-rectangular orthocentric $d$-simplices which are equiradial,
and in the proof of Theorem 4.13 which proves that orthocentric
simplices with $\cI = \cC$ are regular.

\textbf{Theorem 4.8:} \emph{Let $\cH$ be the orthocenter of a
non-rectangular orthocentric $d$-simplex $S = [A_1, \dots,
A_{d+1}]$, and let $\cH'$ be the orthocenter of the face $S' =
[A_1, \dots, A_{k+1}]$ of $S$. Let $a_1, \dots, a_{d+1}$ be the
barycentric coordinates of $\cH$ with respect to $S$, and $a'_1,
\dots, a'_{k+1}$ be the barycentric coordinates of $\cH'$ with
respect to $S'$. Let $\sigma (S)$ be the obtuseness of $S$ as
defined in $($8$)$. Then $a_1 + \dots + a_{k+1} \not= 0$, and}
\begin{equation}
a'_j = \frac{a_j}{a_1 + \cdots + a_{k+1}}\,, \,\,\, \sigma (S') =
\frac{\sigma (S)}{a_1 + \cdots + a_{k+1}}\,.
\end{equation}

\textbf{Proof:} It is clearly sufficient to prove our theorem for $k = d-1$. Also, we may assume
that $\cH$ is the origin $O$. Then it follows from (14) in the proof of Theorem 3.3 that
\[
a'_i = \frac{a_i}{a_1 + \cdots + a_d}\,,
\]
as desired.

To compute $\sigma (S')$, we set $c= \sigma (S)$ and $s = a_1 + \dots + a_d \,(=1-a_{d+1})$. Take
$i$ and $j$ such that $i,j$, and $d+1$ are pairwise distinct. Using (14), we see that
\[
\begin{array}{lll}
A_i \cdot A_j & = & c\\
\cH' \cdot A_i & = & {\displaystyle \frac{-a_{d+1}}{s}} (A_{d+1} \cdot A_i) =
{\displaystyle \frac{-ca_{d+1}}{s}}\\
\|\cH'\|^2 & = & \left({\displaystyle \frac{-a_{d+1}}{s}}\right)^2
\|A_{d+1}\|^2 = \left({\displaystyle \frac{-a_{d+1}}{s}}\right)^2
\left({\displaystyle \frac{-cs}{a_{d+1}}}\right) - {\displaystyle \frac{-ca_{d+1}}{s}}\,,
\end{array}
\]
and therefore
\[
\begin{array}{lll}
\sigma (S') & = & (\cH' - A_i) \cdot (\cH' - A_j)\\
& = & \left({\displaystyle \frac{c}{s}}\right) (-a_{d+1} + 2a_{d+1} + 1 - a_{d+1}) =
{\displaystyle \frac{c}{s}}\,,
\end{array}
\]
as claimed. \hfill $\Box$

\textbf{Theorem 4.9:} \emph{Let $S = [A_1, \dots, A_{d+1}]$ be a
non-rectangular orthocentric $d$-simplex, and let $a_1, \dots,
a_{d+1}$ be the barycentric coordinates of its orthocenter $\cH$.
Let $c = \sigma (S)$ be the obtuseness of $S$ defined in $($8$)$.
Let $\cC$ be the circumcenter and $R$ the circumradius of $S$.
Then}
\begin{eqnarray}
\cC + {\displaystyle \frac{d-1}{2}} \cdot \cH &=& {\displaystyle \frac{1}{2}} (A_1 + \dots + A_{d+1})\,,\\
{\displaystyle \frac{4 R^2}{c}} & = & (d-1)^2 - \sum \limits^{d+1}_{i=1} {\displaystyle \frac{1}{a_i}}\,.
\end{eqnarray}
\emph{Consequently, $\cC$ is interior if and only if $a_i < 1 /
(d-1)$. Also, if $F= [A_1, \dots, A_{k+1}]$ is a face of $S$, and
if $s = s(F) = a_1 + \dots + a_{k+1}$, then the circumradius $R_F$
of $F$ is given by}
\begin{equation}
\frac{4R^2_F}{c} = \frac{(k-1)^2}{s} - \left(\frac{1}{a_1} + \dots + \frac{1}{a_{k+1}}\right)\,.
\end{equation}

\textbf{Proof:} For simplicity, we assume $\cH = O$ (by replacing $A_i$
and $\cC$ by $A_i - \cH$ and $\cC - \cH$, respectively),
and we let $P = (A_1 + \dots + A_{d+1})/2$.
We use the facts that $A_i \cdot A_j = c$ for all $i \not= j$, and that $\|A_i\|^2 = c(a_i-1)/a_i$.
Then $2 (P-A_{d+1}) = A_1 + \dots + A_d - A_{d+1}$ and
\[
\begin{array}{lll}
4 \|P - A_{d+1}\|^2 & = & \|A_1\|^2 + \dots + \| A_{d+1}\|^2 + (d^2 - 3d)c\\
&&\\
&=& {\displaystyle \frac{c(a_1-1)}{a_1}} + \dots + {\displaystyle \frac{c(a_{d+1}-1)}{a_{d+1}}}
+ (d^2 - 3d)c\\
&&\\
&=& -c \left({\displaystyle \frac{1}{a_1}} + \dots + {\displaystyle \frac{1}{a_{d+1}}}\right)
+ (d+1) c + (d^2 - 3d)c\\
&&\\
&=& - c \left({\displaystyle \frac{1}{a_1}} + \dots + {\displaystyle \frac{1}{a_{d+1}}}\right)
+ c (d-1)^2\,.
\end{array}
\]
Similarly for $4\|P-A_i\|^2$ for every $i$. Therefore $P$ is
equidistant from the vertices of $S$, showing that $P$ is the
circumcenter, and proving (28) and (29).

Letting
\[
b_i = \frac{1}{2} + \frac{(1-d)a_i}{2}\,,
\]
we see that
\[
\cC = b_1 A_1 + \dots + b_{d+1} A_{d+1} \,\, \mbox{ and } \,\, b_1 + \dots + b_{d+1} = 1\,.
\]
Therefore $b_1, \dots, b_{d+1}$ are the barycentric coordinates of
$\cC$. Thus $\cC$ is interior if and only if $b_i > 0$, i.e., if
and only if $a_i < 1 / (d-1)$, as claimed.

Finally, the statement pertaining to $R_F$ follows from (29) using
Theorem 4.7. \hfill $\Box$

\subsection{Equiradial orthocentric simplices}

In this subsection we complete the characterization of equiradial
orthocentric simplices, which is given in Theorem 4.9.

We remind the reader that the Gram matrix of a $d$-dimensional
simplex $S=[A_1, \dots, A_{d+1}]$ is the $(d+1) \times
(d+1)$-matrix whose entries $G_{ij}$ are defined by $G_{ij} = A_i
\cdot A_j$, see \cite{H-J}, \cite{L-T}, and our Subsection 1.3.

\textbf{Lemma 4.10:} \emph{Let $G$ be the $(d+1) \times
(d+1)$--matrix whose entries $G_{ij}$ are given by}
\begin{equation}
 G_{ij} = \left\{ \begin{array}{ccc}
 1 & \mbox { if } & i \neq j,\\
 1+x_i & \mbox { if } & i = j. \end{array}\right.
\end{equation}
\emph{where the $x_i$'s are non--zero numbers}.
\begin{itemize}
\item[(a)] \emph{The characteristic polynomial of $G$ is given by}
\begin{equation}
 (\lambda - x_1) \cdots (\lambda-x_{d+1}) \left( 1+\frac{-1}{\lambda-x_1} + \cdots + \frac{-1}{\lambda-x_{d+1}} \right).
\end{equation}
\item[(b)] \emph{Let $x$ and $y$ be distinct non--zero real numbers, and suppose that}
\begin{equation}
x_i = \left\{ \begin{array}{ccc}
 x & \mbox { \emph{if} } & i \le m,\\
 y & \mbox { \emph{if} } & i>m .
\end{array}\right.
\end{equation}

\begin{itemize}
\item[(i)] \emph{If $m=0$, then $cG$ for some $c$ is the Gram
matrix of a $d$--simplex $S$ if and only if $y=-d-1$. In this
case, $S$ is regular}.
\item[(ii)] \emph{If $m=1$, then $cG$ for some $c$ is the Gram
matrix of a $d$--simplex $S$ if and only if}
\begin{equation}
 (i) \,\, xy + dx+my=0 \quad \mbox { \emph{and} } \quad (ii) \,\, x< 0.
\end{equation}
\emph{In this case, $S$ is a kite having eccentricity $(x+y)/(2yx)$}.
\item[(iii)] \emph{If $2 \le m \le d-1$, then $cG$ for some $c$ is
the Gram matrix of a $d$--simplex $S$ if and only if}
\begin{equation}
 (i) \,\, xy + (d+1-m) x+my = 0 \quad \mbox { \emph{and} } \quad (ii) \,\, x+m<0.
\end{equation}
\end{itemize}
\end{itemize}
\textbf{Proof.} Statement (a) follows immediately from Lemma 3.6,
since the characteristic polynomial of $G$ equals
\begin{eqnarray*}
 &  & \det (\lambda I - G)\\[0.3cm]
 & = & \det (\bJ(-1, \cdots, -1; \lambda-x_1, \ldots, \lambda-x_{d+1}))\\[0.3cm]
 & = & (\lambda - x_1) \cdots (\lambda - x_{d+1}) \left( 1 + \frac{-1}{\lambda-x_1} + \cdots + \frac{-1}{\lambda-x_{d+1}} \right) .
\end{eqnarray*}
To prove (b), we use (32) and the fact that $\pm G$ is the Gram matrix of a $d$--simplex if and only if one of its eigenvalues is $0$ and the others all have the same sign.\\
(i) If $m=0$, then the characteristic polynomial of $G$ is
\[
 (\lambda-y)^{d+1} \left( 1 - \frac{d+1}{\lambda - y} \right) = (\lambda - y)^d (\lambda-y-(d+1)),
\]
and the eigenvalues are $y$ and $y+d+1$. Also, $y\neq 0$.
Therefore, $\pm G$ is the Gram matrix of a $d$-simplex $S$ that is
necessarily regular if and only if $y+d+1=0$.

(ii) If $m=1$, then the characteristic polynomial of $G$ is
\begin{eqnarray*}
 f(\lambda) & = & (\lambda - x)(\lambda - y)^d \left( 1-\frac{1}{\lambda - x}- \frac{d}{\lambda - y} \right)\\[0.3cm]
 & = & (\lambda - y)^{d-1} (\lambda^2 - (x+y+d+1)\lambda + y + xy + dx),
\end{eqnarray*}
and one of the eigenvalues is $0$ if and only if
\begin{equation}
 y + xy + dx = 0 .
\end{equation}
Assuming (36), we see that $x+1 \neq 0$ (since $dx \neq 0$), and
that the remaining eigenvalues of $G$ are
\[
 y = \frac{-dx}{x+1} \,\, \mbox { and } \,\, x + y + d+1 =x+\frac{-dx}{x+1} +d+1 = \frac{(x+1)^2+d}{x+1} .
\]
These have the same sign if and only if $x+1$ and $-x(x+1)$ do,
which happens if and only if $x<0$.

One can recover the squares of the edge-length of a simplex from its Gram
matrix by the formula $\|A_i-A_j\|^2 = G_{ii} + G_{jj} - 2 G_{ij}$.
A $d$-simplex $S$ whose Gram matrix is $G$ is thus a $d$-kite whose base is a regular
$(d-1)$-simplex of side-length $(1+y) + (1+y) - 2 = 2y$ and whose remaining edges have
edge-lengths $(1+x) + (1+y) - 2 = (x+y)$. Its eccentricity is therefore $(x+y)/(2y)$.

(iii) If $2 \le m \le d-1$, then, setting $n = d+1-m$, we have $m, n \ge 2$. The
characteristic polynomial of $G$ is
\[
\begin{array}{lll}
f(\lambda) & = & (\lambda-x)^m (\lambda - y)^n \left(1 - {\displaystyle \frac{m}{\lambda - x} -
\frac{n}{\lambda - y}} \right)\\
&&\\
& = & (\lambda-x)^{m-1} (\lambda-y)^{n-1} (\lambda^2 - (x+y+d+1) \lambda + xy + my + nx)\,,
\end{array}
\]
and one of the eigenvalues is $0$ if and only if
\begin{equation}
my + xy + nx = 0\,.
\end{equation}
Assuming (37), we see that $x + m \not= 0$ (since $nx \not= 0$),
and that the remaining eigenvalues of $G$ are
\[
x,\, y = \frac{-nx}{x+m}\,, \mbox{ and } \, x+y+d+1 =
\frac{(x+m)^2 + mn}{x+m}\,.
\]
These have the same sign if and only if $x, -x (x+m)$, and $x+m$
do, which happens if and only if $x+m < 0$. \hfill $\Box$

\subsection{Complete classification of equiradial orthocentric
simplices}

Here we complete the study of orthocentric simplices in which the
circumcenter coincides with the incenter.

\textbf{Theorem 4.11 a:} \emph{Let $S = [A_1, \dots, A_{d+1}]$ be a
non-regular, non-rectangular, equiradial orthocentric $d$-simplex,
and let $a_1, \dots, a_{d+1}$ be the barycentric coordinates of
its orthocenter with respect to $S$. Then there are two
possibilities \newline\newline 1. $S$ is a kite of eccentricity
$\sqrt{(d-2)/d}$, in which case $d\ge 5$, or
\newline 2. There exists an $m$ with $2 \le m \le d-1$ such that, after
relabelling vertices,} \begin{equation} a_1 = \dots = a_m\,, \,\,
a_{m+1} = \dots = a_{d+1}\ \end{equation}\emph{\ where\ }
\begin{equation}m(d+1-m) \le \left(\frac{d^2 -
3d+4}{2(d-2)}\right)^2 , \emph{which implies } d\ge 9.
\end{equation}

\textbf{Proof:} We may assume $d>3$ since equiradial orthocentric
simplices of lower dimensions are known to be regular. Let $S =
[A_1, \dots, A_{d+1}]$ and suppose that the orthocenter is $O$,
and that $a_1, \dots, a_{d+1}$ are the barycenric coordinates of
$O$ with respect to $S$. Suppose that $S$ is equiradial. Then it
follows from (30) in Theorem 4.9 that \[ \frac{(d-2)^2}{1-a_i} +
\frac{1}{a_i} = \frac{(d-2)^2}{1-a_j} + \frac{1}{a_j} \] for every
$i$ and $j$. Therefore \[ \frac{1 + (d-3) (d-1)a_i}{a_i (1-a_i)} =
\frac{1 + (d-3) (d-1)a_j}{a_j (1-a_j)}\,, \] which simplifies into
\[ (a_i - a_j) ((d-3) (d-1) a_i a_j + a_i + a_j - 1) = 0\,. \] If
$(d-3) (d-1) a_i a_j + a_i + a_j = 1$ and $(d-3) (d-1) a_i a_k +
a_i + a_k = 1$ then, multiplying the first by $a_k$ and the second
by $a_j$ and subtracting, we obtain $(a_i - 1) (a_k - a_j) = 0$,
and therefore $a_k = a_j$, since no $a_i$ can be $1$. Therefore
the $a_i$'s can take at most two different values $a$ and $b$ that
satisfy
\begin{equation}
(d-3) (d-1) ab + a + b -1 = 0\,.
\end{equation}
Since $S$ is not regular, it follows that the $a_i$'s are not all
equal, and therefore we may assume that
\[ a_1 = \dots = a_m = a\,, \,\,\, a_{m+1} = \dots =
a_{d+1} = b
\]
for some $m$ with $1 \le m \le d$, and with $a$ and
$b$ satisfying (40). Letting $n = d+1-m, \,\, x = -1/a$ and $y =
-1/b$, and using $ma + nb = 1$ and (40), we obtain (37) and
\begin{equation}
xy + x + y = (d-3) (d-1)\,.
\end{equation}
Also, using $A_i \cdot A_j = c$ for $i \not= j$, and
$\|A_i\|^2 = c(a_i - 1)/a_i$, we see that the Gram matrix of $S$
is of the type described in Lemma 4.10 (b), where $x$ and $y$
satisfy (41).

Again we use Lemma 4.10 (b). If $m = 1$, then it follows from (34)
and (41) that $x = 3-d$ and $x < 0$. Substituting $m=1, n=d$, and
$x=3-d$ in (37), we see that $d \ge 5$ and that
\[
y = \frac{d(d-3)}{d-4}\,.
\]
This corresponds to the kite whose eccentricity $\epsilon$ is
given by
\[ \epsilon^2 = \frac{x+y}{2y} = \frac{(d-2)^2}{d^2}\,,
\]
in conformance with Theorem 4.7. Thus we are left with the case
\[
2 \le m \le d-1\,,
\]
and
\begin{equation} x+m < 0\,, \,\, xy +
my + nx = 0\,, \,\, x+y+xy = (d-3) (d-1)\,.
\end{equation}
The last two of these can be rewritten as
\[
\begin{array}{rrl} (x+m) (y+n) & = & mn\,,\\ &&\\ (x+m) (1-n) + (y+n)
(1-m) & = & d^2 - 3d + 4 - 2mn\,.
\end{array}
\]
Setting
\begin{equation}
\xi = (x+m)(1-n), \eta = (y+n) (1-m)\,,
\end{equation}
we see that (42) can be rewritten as
\[ \xi > 0, \,
\eta > 0, \, \xi + \eta = d^2 - 3d+4 - 2mn, \xi \eta = mn
(mn-d)\,.
\]

Now
\[
0\le (\xi - \eta)^2 = (\xi + \eta)^2 - 4 \xi \eta = (d^2 - 3d + 4
- 2mn)^2 - 4mn (mn-d)\,.
\]
This simplifies into
\begin{equation}
(d^2 - 3d+4)^2 - 4mn (d^2 - 3d+4) + 4mnd \ge 0
\end{equation}
or, equivalently,
\begin{equation}
mn \le \frac{(d^2 - 3d+4)^2}{4(d-2)^2}\,,
\end{equation}
as desired. It remains to prove that $d \ge 9$. Since $mn \ge 2
(d-1)$, it follows that
\begin{equation}
2 (d-1) \le \frac{(d^2 - 3d + 4)^2}{4(d-2)^2}\,,
\end{equation}
and therefore
\[
0 \le d^4 - 14d^3 + 57d^2 - 88d + 48 = (d-3) \left(d(d-3) (d-8) -
16\right) \,. \] Since $d>3$ by assumption, this happens if and
only if $d \ge 9$, as claimed.\hfill $\Box$

The kites have already been completely analyzed in Subsection 4.2.
In the other case we have the following converse.

\textbf{Theorem 4.11 b:}  \emph{If $2 \le m \le d-1$, and if $m$
and $d$ satisfy $(39)$, then there exist exactly two non-similar
non-rectangular orthocentric equiradial $d$-simplices whose
orthocenter's barycentric coordinates $a_1, \dots, a_{d+1}$
satisfy $(38)$. This happens for any given value of $m$ if $d$ is
large enough; in particular when $m =2$ and $d \ge 9$}.

\textbf{Proof:} Suppose that $2 \le m \le d-1$ and that (39)
holds. Note first that equality cannot take place in (39), since
$(d^2-3d+4)/(2(d-2))$  is not an integer. This can be seen by
taking the cases $d \equiv 1 ~(\mbox{mod} ~2)$, $d \equiv 0
~(\mbox{mod}~ 4)$, and $d \equiv 2 ~(\mbox{mod}~ 4)$. Then it
follows immediately that  the discriminant of
\begin{eqnarray*}
Q(Z):= Z^2 - (d^2-3d+4-2mn)Z+mn(mn-d)
\end{eqnarray*}
is strictly positive. Therefore $Q(Z)$ has two distinct real zeros
which are necessarily positive since their sum $d^2-3d+4-2mn$ and
product $mn(mn-d)$ are. Letting $X$ and $Y$ be the zeros of
$Q(Z)$, we find $x$ and $y$ by solving  the system
\begin{equation}
X = (x+m) (1-n), \, Y= (y+m) (1-n)\,
\end{equation}
or the system
\begin{equation}
Y = (x+m) (1-n), \, X= (y+m) (1-n)\,.
\end{equation}
Since $X$ and $Y$ are distinct, these systems give rise to two
distinct pairs $(x,y)$. The Gram matrices corresponding to these values give
rise to the desired simplices, as in the proof of Lemma 4.10.
This completes the proof. \hfill $\Box$

\textbf{Remark 4.12:} Consideration of the corresponding Gram matrices shows
that these orthocentric equiradial $d$-simplices may be thought of as generalized
kites. They may be described as the join of a regular $(m-1)$-simplex of
edge-length $a$ with a regular $(d-m)$-simplex of edge-length $b$ such that all
intervening edges have edge-length $c$, for suitable values of the parameter $d, m,
a, b, c$.

\subsection{Orthocentric simplices with circumcenter = incenter
are regular}

Here we complete the study of orthocentric simplices in which the
circumcenter coincides with the incenter.

\textbf{Theorem 4.13:} \emph{If $S$ is an orthocentric $d$-simplex
in which the circumcenter and the incenter coincide, then $S$ is
regular}.

\textbf{Proof:} We shall treat the simpler case when $S$ is
rectangular in Section 5. So we suppose that $S$ is a
non-rectangular orthocentric $d$-simplex. Let $a_1, \dots,
a_{d+1}$ be the barycentric coordinates of its orthocenter with
respect to $S$. If the incenter and the circumcenter of $S$
coincide, and if $S$ is not regular, then, being necessarily
equiradial, $S$ is of one of the two types in Theorem 4.11. The
first is a $d$-kite of eccentricity $\sqrt{(d-2)/d}$, and
therefore has an exterior circumcenter by Theorem 4.7. The other
type satisfies
\[
a_1 = \dots = a_m = a\,, \,\,\, a_{m+1} = \dots = a_{d+1} = b\,,
\]
where $x = -1/a$ and $y = -1/b$ are such that
\begin{equation}
x + m < 0\,, \,\,\, y + n < 0\,, \,\,\, xy + x + y = (d-3) (d-1)\,.
\end{equation}
In particular,  $x$ and $y$ are negative and
\[
\begin{array}{lll}
\mbox{the circumcenter of } S \mbox{ is interior } & \Longleftrightarrow & a, b < {\displaystyle \frac{1}{d-1}},
\mbox{by Theorem 3.10,}\\
&&\\
& \Longleftrightarrow & x,y < 1 - d\\
&&\\
& \Longleftrightarrow & x + 1, y + 1 < 2-d\\
&&\\
&\Longrightarrow & (x+1) (y+1) > (d-2)^2\\
&&\\
& \Longrightarrow & xy + x + y > (d-1) (d-3)\,,
\end{array}
\]
contradicting (49). Therefore here again the circumcenter is not
interior, and cannot coincide with the incenter. \hfill $\Box$

\textbf{Remark 4.14:} We have seen in \cite{E-H-M}, Theorem 2.2,
that there are non-equifacetal tetrahedra whose facets have equal
inradii. This prompts the still open question whether there exist
non-regular orthocentric tetrahedra (or higher dimensional
simplices) with this property.

\section{Rectangular simplices}

In the previous section, we have coordinatized a non-rectangular
orthocentric $d$-simplex $S = [A_1, \dots, A_{d+1}]$ by the
barycentric coordinates $a_1, \dots, a_{d+1}$ of its orthocenter,
together with the obtuseness $\sigma (S)$ defined in (8). It was
noted that $\sigma (S) = 0$ if and only if $S$ is rectangular, and
it is easy to see that if $S'$ is similar to $S$ with a similarity
factor $p$, then $\sigma (S') = p^2 \sigma (S)$. Thus $|\sigma
(S)|$ is not relevant as far as the shape of $S$ is concerned, and
$S$ can be scaled so that $\sigma (S) = 0$ or $\sigma (S) = \pm
1$. In view of (10) it is obvious that $\sigma (S) < 0$ if and
only if $0 < a_i < 1$ for all $i$, i.e., if and only if the
orthocenter of $S$ is an interior point. Note that, for a
triangle, this is equivalent to $S$ being acute-angled.

Rectangular $d$-simplices are characterized as those orthocentric
$d$-simplices $S$ with $\sigma (S) = 0$. However, the barycentric
coordinates of the orthocenter of such a $d$-simplex $S$ carry no
information about $S$ (except for locating the vertex at which the
orthocenter occurs), and therefore cannot serve to parametrize
$S$. On the other hand, the lengths $b_1, \dots, b_d$ of the legs
of a rectangular simplex do carry all the essential information
about $S$. Here a \emph{leg} is an edge emanating from the vertex
at which $S$ is rectangular. If a $d$-simplex $ S = [A_1, \dots,
A_{d+1}]$ is rectangular, say at $A_{d+1}$, then one can place it
in $\bE^d$ in such a way that $A_{d+1}$ lies at the origin and
such that the legs $A_{d+1} A_i$ lie on the positive coordinate
axes. Then the $i$-th cartesian coordinate of $A_i$ is $b_i$ and
the other coordinates are $0$, and
\begin{equation}
 \Vert A_{d+1} - A_i \Vert = b_i .
\end{equation}
The volume, the circumradius, the inradius, and other elements of
$S$ can be easily computed in terms of the $b_i$'s, as illustrated
in Theorem 5.1 below. These formulas will then apply to, but only
to, the rectangular faces of $S$, i.e., to those faces of $S$ that
have $A_{d+1}$ as a vertex. To understand the remaining faces,
note that they are necessarily faces of the \textit{hypotenuse
facet} $[A_1, \cdots, A_d]$, and thus it is sufficient and
important to understand this facet. It follows from (57) below
that this facet, which is necessarily orthocentric, cannot be
rectangular, and therefore it yields to the results of the
preceding section. The natural question that arises is whether
every non--rectangular orthocentric $(d-1)$-simplex occurs as a
facet (necessarily the hypotenuse facet) of a rectangular
$d$-simplex. Theorem 5.3 below provides a satisfactory answer.

\textbf{Theorem 5.1} \emph{Let $S = [A_1, \cdots, A_{d+1}]$ be a
$d$-simplex that is rectangular at $A_{d+1}$, and let $b_1,
\ldots, b_d$ be the lengths of its legs $A_{d+1} A_1, \ldots,
A_{d+1}A_d$, respectively. Let the volume, the circumradius, and
the inradius of $S$ be denoted by $V, R$, and $r$, respectively.
For each $i$, let $V_i$ be the $(d-1)$-volume of the $i$-th
facet of $S$, and let $h$ be the altitude of $S$ to the
$(d+1)$-th facet $[A_1, \cdots, A_d]$. Then}
\begin{eqnarray}
V & = & \frac{b_1 \cdots b_d}{d!}\,, \\[0.3cm]
V_{d+1} & = & \frac{b_1 \cdots b_d}{(d-1)!} \, \sqrt{\frac{1}{b^2_1} + \cdots + \frac{1}{b^2_d}}\,,\\[0.3cm]
h & = & \left(  \sqrt{\frac{1}{b^2_1} + \cdots + \frac{1}{b^2_d}} \right)^{-1}\,, \\[0.3cm]
r & = & \left( \frac{1}{b_1} + \cdots + \frac{1}{b_d} +
\sqrt{\frac{1}{b^2_1}
+ \cdots + \frac{1}{b^2_d}} \right)^{-1} \,,\\[0.3cm]
R^2 & = & \frac{b^2_1 + \cdots + b^2_d}{4}\,.
\end{eqnarray}
\emph{Also, if $A_{d+1} = O$, then the circumcenter $\cC$ of $S$ and the orthocenter $\cB$
of the facet $[A_1, \cdots, A_d]$ are given by}
\begin{eqnarray}
\cC & = & \frac{A_1 + \cdots + A_d}{2}\,,\\[0.3cm]
 \cB & = & \left( \frac{1}{b^2_1} + \cdots + \frac{1}{b^2_d}\right)^{-1}
 \left( \frac{1}{b^2_1}A_1  + \cdots + \frac{1}{b^2_d} A_d \right)\,.
\end{eqnarray}

\textbf{Proof:} The equation (51) is obvious. Using (51) and the
$d$-dimensional Pythagoras' Theorem, we obtain
\[
\begin{array}{lll}
V^2_{d+1} & = & \left({\displaystyle \frac{b_1 \cdots
b_d}{(d-1)!\,b_d}}\right)^2 + \cdots +
\left( {\displaystyle \frac{b_1 \cdots b_d}{(d-1)!\,b_d}}\right)^2\\
&&\\
& = & \left({\displaystyle \frac{b_1 \cdots b_d}{(d-1)!}}\right)^2
\left({\displaystyle \frac{1}{b^2_1} + \dots + \frac{1}{b^2_d}}
\right)\,.
\end{array}
\]
This proves (52). Then we use $dV = V_{d+1} h$ to get (53). For
(54), we use the preceding formulas and the fact that $dV = r (V_1
+ V_2 + \dots + V_{d+1})$. To prove (56), we use $\|\cC\|^2 =
\|\cC - A_i\|^2$ to conclude that $2\cC \cdot A_i = \|A_i\|^2 =
b^2_i$, and therefore the $i$-th coordinate of $\cC$ is $b_i/2$.
Then we use (56) and Pythagoras' Theorem to obtain (55). Finally,
(57) follows from the fact that $\cB$ is the projection of
$A_{d+1}$ on the facet $[A_1, \dots, A_d]$ and thus $\cB \cdot
(A_i - A_j) = 0$ for $1 \le i < j \le d$. \hfill $\Box$

\textbf{Theorem 5.2:} \emph{If $S$ is rectangular, then its four
classical centers $\cG, \cC, \cI$, and $\cH$ are pairwise
distinct}.

\textbf{Proof:} Let $S = [A_1, \dots, A_{d+1}]$ be rectangular at
$A_{d+1}$. In view of (55) the circumcenter lies on the hypotenuse
facet of a rectangular $d$-simplex if and only if $d=2$. Also, the
barycentric coordinates of $\cC$ with respect to $A_1, \dots,
A_{d+1}$ are $(1/2, \dots, 1/2, (2-d)/2)$, and therefore $\cC$ is
never interior (since $(2-d)/2 \le 0$ for all $d \ge 2$).
Therefore, of the orthocenter, the circumcenter, the incenter, and
the centroid of $S$, the only two that can possibly coincide are
the last two. This cannot happen either. In fact, it is clear that
$\cI = (r, \dots, r)$ and $\cG = (b_1, \dots, b_{d+1})/(d+1)$, and
therefore the equality $\cI = \cG$ would imply that $b_i = (d+1)
r$ for all $i$. Using (53), we arrive at the contradiction
\begin{equation}
d+1 = d + \sqrt{d} \,.
\end{equation}
Therefore no two the four classical centers can coincide. \hfill $\Box$

\textbf{Theorem 5.3:} \emph{Let $T = [A_1, \dots, A_d]$ be an
orthocentric $(d-1)$-simplex. Then the following conditions are
equivalent.}
\begin{enumerate}
\item[(a)] \emph{There exists a $d$-simplex $S = [A_1, \dots, A_d, A_{d+1}]$ that is
rectangular at $A_{d+1}$.}
\item[(b)] \emph{The orthocenter of ~$T$ is interior}.
\item[(c)] \emph{$\sigma (T) < 0$.}
\end{enumerate}

\textbf{Proof:} We already know that (b) and (c) are equivalent by
Theorem 3.8. Let $t_1, \dots, t_d$ be the barycentric coordinates
of the orthocenter $\cB$ of $T$. If (a) holds, then (57) implies
that the $t_i$'s are all positive, and therefore $\cB$ is
interior. Thus $0 < t_i < 1$ for all $i$. In view of (10), this
implies that $\sigma (T) < 0$. Conversely, if $c = \sigma (T) <
0$, then $0 < t_i < 1$ for all $i$. Let $b_1, \dots, b_d$ be the
positive numbers defined by
\begin{equation}
b^2_i t_i = -c\,,
\end{equation}
and let $S = [P_1, \dots, P_{d+1}]$ be the $d$-simplex that is rectangular at $P_{d+1}$ and whose
leg-lengths are $b_1, \dots, b_d$. Then for $1 \le i < j \le d$ we have
\[
\begin{array}{lll}
\|A_i - A_j\|^2 & = & -c \left({\displaystyle \frac{1}{t_i} +
\frac{1}{t_j}}\right)\,\,\,\,
\mbox{ by (11)}\\
&&\\
& = & b^2_i + b^2_j\,\,\,\, \mbox{ by (59)}\\
&&\\
& =& \|P_i - P_j\|^2\,\,\,\, \mbox{ by Pythagoras' Theorem}.
\end{array}
\]
Therefore $T$ is congruent to the facet $[P_1, \dots, P_d]$ of $S$. This proves (a).
\hfill $\Box$

\end{document}